\documentclass{article}
\usepackage{bm,latexsym,amssymb,amsthm,amsmath}
\usepackage{titlesec}
\usepackage{array,contour,color,graphicx,booktabs,threeparttable,tabularx}
\usepackage[dvipsnames]{xcolor}
\usepackage{fancyhdr,makeidx}
\begin{document}
\title{Nonstandard Vector Space with a Metric and Its Topological Structure}
\author{\vspace{2mm}Hsien-Chung Wu
\thanks{e-mail:\vspace{5mm}
hcwu@nknucc.nknu.edu.tw}\\
Department of Mathematics, \\
National Kaohsiung Normal University, Kaohsiung 802, Taiwan}
\date{}
\maketitle
\newtheorem{Thm}{Theorem}[section]
\newtheorem{Def}{Definition}[section]
\newtheorem{Lem}{Lemma}[section]
\newtheorem{Pro}{Proposition}[section]
\newtheorem{Rem}{Remark}[section]
\newtheorem{Cor}{Corollary}[section]
\newtheorem{Ex}{Example}[section]

\newenvironment{Proof}{\noindent {\bf Proof}.}
{\hspace{5mm}\rule[-0.8mm]{2.5mm}{2mm}\medskip\par}

\begin{abstract}
In this paper, we introduce the nonstandard vector space
in which the concept of additive inverse element will not be taken into account.
We also consider a metric defined on this nonstandard vector space.
Under these settings, the conventional intuitive properties for the open 
and closed balls will not hold true. Therefore, four kinds of open and closed 
sets are proposed. Furthermore, the topologies generated by these different 
concepts of open and closed sets are investigated.

\vspace{3 mm}


\vspace{3 mm}

\end{abstract}

\section{Introduction}

It is well-known that the topic of topological vector space is based on the 
vector space by referring to the monographs \cite{ada,kha,sch,per,won}. 
We can see that the set of all closed intervals cannot form 
a real vector space. The main reason is that there will be no additive 
inverse element for each closed interval in $\mathbb{R}$. Under this inspiration, 
we are going to propose the concept of nonstandard vector space 
which will not own the concept of additive inverse element and 
will be weaker than the (conventional) vector space in the sense of axioms. 

Let $X$ be an universal set. We can define a (conventional) metric $d$ on $X$, 
which satisfies the following conditions:
\begin{itemize}
\item for any $x,y\in X$, $d(x,y)=0$ if and only if $x=y$;

\item for any $x,y\in X$, $d(x,y)=d(y,x)$;

\item for any $x,y,z\in X$, $d(x,y)\leq d(x,z)+d(z,y)$.
\end{itemize}
If $X$ is taken as a vector space, then we can define a norm 
$\parallel\cdot\parallel$ on $X$, which satisfies the following 
conditions:
\begin{itemize}
\item $\parallel x\parallel$ if and only if $x=\theta$;

\item $\parallel\alpha x\parallel =|\alpha |\cdot\parallel x\parallel$
for any $x\in X$ and $\alpha\in\mathbb{F}$;

\item for any $x,y\in X$, $\parallel x+y\parallel\leq
\parallel x\parallel +\parallel y\parallel$.
\end{itemize}
It is well-known that the normed space $(X,\parallel\cdot\parallel )$
is also a metric space $(X,d)$ with the metric $d$ defined by 
$d(x,y)=\parallel x-y\parallel$ such that $d$ is translation-invariant
and homogeneous, i.e., $d(x+a,y+a)=d(x,y)$ and $d(\alpha x,\alpha y)=
|\alpha |\cdot d(x,y)$.

Conversely, we assume that the metric space $(X,d)$ is also a vector space.
If $d$ is not translation-invariant and homogeneous, then we cannot induced
a normed space from $(X,d)$ by defining a suitable norm based on the metric $d$.
Of course, if $d$ is translation-invariant and homogeneous and we
define a nonnegative function $p(x)=d(x,\theta )$, then we can show that
$p$ is indeed a norm on $X$. Therefore, when $X$ is taken as 
a vector space and the metric $d$ is not translation-invariant and 
homogeneous, the metric space $(X,d)$ is weaker than the normed space 
$(X,\parallel\cdot\parallel )$ in the sense that the properties which 
hold true in $(X,\parallel\cdot\parallel )$ are also hold true in 
$(X,d)$ with $d(x,y)=\parallel x-y\parallel$, and the converse is not true.
When $X$ is taken as the nonstandard vector space over 
$\mathbb{F}$, we are going to study the so-called nonstandard metric space 
in which $X$ is endowed with a metric $d$. Therefore, in the nonstandard 
metric space, we can perform the vector addition and scalar multiplication.
In the conventional metric space, the vector addition and scalar multiplication
are not allowed, unless the universal set $X$ is taken as the vector space.
We have to mention that the mathematical structures of nonstandard metric 
space that is based on the nonstandard vector space are completely different from 
that of the (conventional) metric space that is based on the (conventional) 
vector space. This is the main purpose of this paper.

In Sections 2 and 3, the concept of nonstandard vector space is proposed, 
where some interesting properties are derived in order to study the 
the topology generated by this kind of space.
In Section 4, we introduce the concept of metric defined on the nonstandard vector 
space defined in Section 2. In Section 5, we provide the non-intuitive properties
for the open and closed balls. In Sections 6 and 7, we propose many different 
concepts of nonstandard open and closed sets. Finally, in Section 8, 
we investigate the topologies generated by these different 
concepts of open and closed sets.

\section{Nonstandard Vector Spaces}

Let $X$ be a universal set and let$\mathbb{F}$ be a scalar field.
We assume that $X$ is endowed with the vector addition 
$x\oplus y$ and scalar multiplication $\alpha x$ for any 
$x,y\in X$ and $\alpha\in\mathbb{F}$. If $X$ is also closed under the vector addition 
and scalar multiplication, then we call
$X$ a {\em universal set over} $\mathbb{F}$. In the conventional 
vector space over $\mathbb{F}$, the additive inverse element of $x$ is 
denoted by $-x$, and it can also be shown that $-x=-1x$. Here, we shall not 
consider the concept of inverse element. However, for convenience, 
we still adopt $-x=-1x$. 

For $x,y\in X$, the {\em substraction} $x\ominus y$ is defined by 
$x\ominus y=x\oplus (-y)$, where $-y$ means the scalar multiplication $(-1)y$.
For any $x\in X$ and $\alpha\in\mathbb{F}$, 
we have to mention that $(-\alpha )x\neq -\alpha x$ and 
$\alpha (-x)\neq -\alpha x$ in general, unless $\alpha (\beta x)=
(\alpha\beta )x$ for any $\alpha ,\beta\in\mathbb{F}$.
Here, this law will not always be assumed to be true.

\begin{Ex}{\label{nhbex1}}
{\em
Let ${\cal I}$ be the set of all closed intervals in $\mathbb{R}$. The vector
addition is given by 
\[[a,b]\oplus [c,d]=[a+c,b+d]\]
and the scalar multiplication is given by
\[k[a,b]=\left\{\begin{array}{ll}
\mbox{$[ka,kb]$} & \mbox{if $k\geq 0$}\\
\mbox{$[kb,ka]$} & \mbox{if $k<0$.}
\end{array}\right .\]
It is easy to see that ${\cal I}$ is not a (conventional) vector space under
the above vector addition and scalar multiplication. The main reason
is that the inverse element does not exist for any non-degenerated closed interval.
}\end{Ex}

Since the universal set $X$ over $\mathbb{F}$ can just own the 
vector addition and scalar multiplication, in general, the universal
set $X$ cannot own the zero element. The set 
\[\Omega =\{x\ominus x:x\in X\}\] 
is called the {\em null set} of $X$. Therefore, the null set can be regarded as 
a kind of zero element of $X$. 

\begin{Ex}{\label{nhbex2}}
{\em
Continued from Example~\ref{nhbex1}, for any $[a,b]\in {\cal I}$, we have 
\[[a,b]\ominus [a,b]=[a,b]\oplus (-[a,b])=[a,b]\oplus [-b,-a]=[a-b,b-a]=[-(b-a),b-a].\]
Therefore, we have the null set $\Omega =\{[-k,k]:k\geq 0\}$.
}\end{Ex}

Now, we are in a position to define the concept of nonstandard vector space.

\begin{Def}{\label{ch2d521}}
{\em
Let $X$ be a universal set over $\mathbb{F}$. 
We say that $X$ is a {\em nonstandard vector space} over $\mathbb{F}$ if 
and only if the following conditions are satisfied:
\begin{description}
\item (i) $1x=x$ for any $x\in X$;

\item (ii) $x=y$ implies $x\oplus z=y\oplus z$ and $\alpha x=\alpha y$
for any $x,y,z\in X$ and $\alpha\in\mathbb{F}$.

\item (iii) the commutative and associative laws for vector 
addition hold true in $X$; that is, $x\oplus y=y\oplus x$ and $(x\oplus y)
\oplus z=x\oplus (y\oplus z)$ for any $x,y,z\in X$.
\end{description}
}\end{Def}

Let $X$ be a universal set over $\mathbb{F}$. 
More laws about the vector addition and scalar multiplication can be defined below.
\begin{itemize}
\item We say that the 
{\em distributive law for vector addition} holds true
in $X$ if $\alpha (x\oplus y)=\alpha x\oplus\alpha y$ for any $x,y\in X$
and $\alpha\in\mathbb{F}$. 

\item We say that the {\em positively distributive law for vector addition} 
holds true in $X$ if $\alpha (x\oplus y)=\alpha x\oplus\alpha y$ 
for any $x,y\in X$ and $\alpha >0$. 

\item We say that the 
{\em associative law for scalar multiplication} holds true in $X$ 
if $\alpha (\beta x)=
(\alpha\beta )x$ for any $x\in X$ and $\alpha ,\beta\in\mathbb{F}$.

\item We say that the 
{\em associative law for positive scalar multiplication} holds true in $X$ 
if $\alpha (\beta x)=
(\alpha\beta )x$ for any $x\in X$ and $\alpha ,\beta >0$.

\item We say that the {\em distributive law for 
scalar addition} holds true in $X$ if $(\alpha +\beta )x=\alpha x\oplus\beta x$ 
for any $x\in X$ and $\alpha ,\beta\in\mathbb{F}$.

\item We say that the {\em distributive law for positive 
scalar addition} holds true in $X$ if $(\alpha +\beta )x=\alpha x\oplus\beta x$ 
for any $x\in X$ and $\alpha ,\beta >0$.

\item We say that the {\em distributive law for negative
scalar addition} holds true in $X$ if $(\alpha +\beta )x=\alpha x\oplus\beta x$ 
for any $x\in X$ and $\alpha ,\beta <0$.
\end{itemize}

We remark that if the distributive law for positive and negative
scalar addition hold true in $X$, then $(\alpha +\beta )x=\alpha x\oplus\beta x$ 
for any $x\in X$ and $\alpha\beta >0$. 

\begin{Ex}{\label{nmsex10}}
{\em 
It is not hard to see that 
the distributive law for scalar addition does not hold true in ${\cal I}$
(the set of all closed intervals). However, 
the distributive law for positive and negative
scalar addition hold true in ${\cal I}$.
}\end{Ex}

Let $X$ be a nonstandard vector space over $\mathbb{F}$
with the null set $\Omega$. We say that $\theta$ is the {\em zero element} of $X$
if and only if $x=\theta\oplus x=x\oplus\theta$ for each $x\in X$.
We write $x\stackrel{\Omega}{=}y$ if one of the following conditions is satisfied:
\begin{enumerate}
\item [(a)] $x=y$;

\item [(b)] there exists $\omega\in\Omega$ such that $x=y\oplus\omega$ or 
$x\oplus\omega =y$;

\item [(c)] there exist $\omega_{1},\omega_{2}\in\Omega$ such that 
$x\oplus\omega_{1}=y\oplus\omega_{2}$.
\end{enumerate}

\begin{Rem}
{\em
Suppose that the nonstandard vector space $X$ also contains the zero 
element $\theta$. Then we can simply say that $x\stackrel{\Omega}{=}y$ 
if and only if there exist $\omega_{1},\omega_{2}\in\Omega$ such that 
$x\oplus\omega_{1}=y\oplus\omega_{2}$ (i.e., only condition (c) is 
satisfied), since conditions (a) and (b) can be rewritten as condition 
(c) by adding $\theta$. We also remark that if we want to discuss some 
properties based on $x\stackrel{\Omega}{=}y$, it suffices to just consider the 
case of condition (c) $x\oplus\omega_{1}=y\oplus\omega_{2}$, even though $X$ does not 
contain the zero elment $\theta$. The reason is that the same arguments 
are still applicable for the cases of condition (a) or (b)
when condition (c) has been shown to be valid.
}\end{Rem}

For any $x,y,z\in X$, we see that $x\oplus z=y\oplus z$ does not imply
$x=y$, i.e., the cancellation law does not hold true.
However, we have the following interesting results.

\begin{Pro}{\label{ch1p*106}}
Let $X$ be a nonstandard vector space over $\mathbb{F}$
with the null set $\Theta$. The following statements hold true.
\begin{enumerate}
\item [{\em (i)}] For any $x,y,z\in X$, if $x\oplus z=y\oplus z$, 
then $x\stackrel{\Omega}{=}y$.

\item [{\em (ii)}] If $x\ominus y\in\Omega$, then $x\stackrel{\Omega}{=}y$. 

\item [{\em (iii)}] Suppose that $\Omega$ is closed under the vector addition.
If $x\stackrel{\Omega}{=}y$, then there exists $\omega\in\Omega$ 
such that $x\ominus y\oplus\omega\in\Omega$.
\end{enumerate}
\end{Pro}
\begin{Proof}
Since $x\oplus z=y\oplus z$, by adding 
$-z$ on both sides, we obtain $x\oplus\omega =y\oplus\omega$, where 
$\omega =z\ominus z\in\Omega$. This proves (i). For result (ii),
we have $x\oplus (-y)=\omega_{1}$ for some $\omega_{1}\in\Omega$. Therefore
$x\oplus (-y)\oplus y=\omega_{1}\oplus y$ by adding $y$ on both sides.
This shows that $x\oplus\omega_{2}=\omega_{1}\oplus y$ for $\omega_{1},
\omega_{2}\in\Omega$, which proves (ii). Now suppose that $x\stackrel{\Omega}{=}y$.
Then $x\oplus\omega_{2}=\omega_{1}\oplus y$ for soem $\omega_{1},
\omega_{2}\in\Omega$. By adding $-y$ on both sides, we obtain 
$x\ominus y\oplus\omega_{2}=\omega_{1}\oplus\omega_{3}\in\Omega$, where 
$\omega_{3}=y\ominus y\in\Omega$. We complete the proof.
\end{Proof}

If $X$ is a nonstandard vector space over $\mathbb{F}$, 
then $-(x\oplus y)=(-x)\oplus (-y)$ does not hold true in general. 
However, we have the following interesting results.

\begin{Pro}{\label{ch1p*102}}
Let $X$ be a nonstandard vector space over $\mathbb{F}$ 
with the null set $\Omega$. 
Suppose that $\Omega$ is closed under the vector addition.
We have $-(x_{1}\oplus\cdots\oplus x_{n})\stackrel{\Omega}{=}
(-x_{1})\oplus\cdots\oplus (-x_{n})$ and 
$-(x\ominus y)\stackrel{\Omega}{=}(-x)\oplus y$.
\end{Pro}
\begin{Proof}
Since $\Omega$ is closed under the vector addition, we have 
$x_{1}\oplus\cdots\oplus x_{n}\oplus (-x_{1})\oplus\cdots
\oplus (-x_{n})=\omega_{1}$ for some $\omega_{1}\in\Omega$. 
By addding $-(x_{1}\oplus\cdots\oplus x_{n})$ on both sides, we obtain 
\[-(x_{1}\oplus\cdots\oplus x_{n})\oplus (x_{1}\oplus\cdots\oplus x_{n})
\oplus (-x_{1})\oplus\cdots\oplus (-x_{n})
=-(x_{1}\oplus\cdots\oplus x_{n})\oplus\omega_{1},\]
which says that $\omega_{2}\oplus (-x_{1})\oplus\cdots\oplus (-x_{n})=
-(x_{1}\oplus\cdots\oplus x_{n})\oplus\omega_{1}$,
where $\omega_{2}=-(x_{1}\oplus\cdots\oplus x_{n})\oplus 
(x_{1}\oplus\cdots\oplus x_{n})\in\Omega$. Therefore, we obtain
$(-x_{1})\oplus\cdots\oplus (-x_{n})\stackrel{\Omega}{=}
-(x_{1}\oplus\cdots\oplus x_{n})$. Similarly, we have
$x\ominus y\oplus (-x)\oplus y=\omega_{3}$ for some $\omega_{3}\in\Omega$. 
Therefore, we have $\omega_{4}\oplus (-x)\oplus y=-(x\ominus y)
\oplus\omega_{3}$, which also means $-(x\ominus y)\stackrel{\Omega}{=}
(-x)\oplus y$. We complete the proof.
\end{Proof}

Now, we are going to introduce the concepts of generalized inverse elements.
Also, the uniqueness of inverse element is, in some sense, different from 
that of conventional vector space.

\begin{Def}
{\em
Let $X$ be a nonstandard vector space over $\mathbb{F}$ with the null set $\Omega$. 
For any $x\in X$, we say that $y$ is the {\em generalized inverse} of $x$ 
if and only if $x\oplus y\in\Theta$.
}\end{Def}

\begin{Pro}{\label{ch2p75}}
Let $X$ be a nonstandard vector space over $\mathbb{F}$
with the null set $\Omega$. For any $x\in X$, if $y$ and $z$ are the 
generalized inverse of $x$, then $y\stackrel{\Omega}{=}z$.
\end{Pro}
\begin{Proof}
Suppose that $y$ and $z$ are the generalized inverse.
Then we have $x\oplus y=\omega_{1}$ and $x\oplus z=\omega_{2}$ 
for some $\omega_{1},\omega_{2}\in\Omega$. By adding $z$ on both sides
of the first equality, we have $x\oplus z\oplus y=z\oplus\omega_{1}$,
i.e., $y\oplus\omega_{2}=z\oplus\omega_{1}$. This shows that 
$y\stackrel{\Omega}{=}z$. We complete the proof.
\end{Proof}

\begin{Cor}
Let $X$ be a nonstandard vector space over $\mathbb{F}$
with the null set $\Omega$. 
For any $x\in X$, if $y$ is the generalized inverse 
of $x$, then $y\stackrel{\Omega}{=}-x=-1x$.
\end{Cor}
\begin{Proof}
By definition, we have $x\ominus x=x\oplus (-x)\in\Theta$. 
This shows that $-x=-1x$ is the generalized inverse of $x$. 
The results follow from Proposition~\ref{ch2p75} immediately.
\end{Proof}

Let $X$ be a nonstandard vector space over $\mathbb{F}$, 
and let $Y$ be a subset of $X$. 
We say that $Y$ is a {\em subspace} of $X$ if and only if $Y$ is closed under 
the vector addition and scalar multiplication, i.e., 
$x\oplus y\in Y$ and $\alpha x\in Y$ for any $x,y\in Y$ 
and $\alpha\in\mathbb{F}$.

\begin{Rem}{\label{ch1r242}}
{\em
Let $Y$ be a subspace of $X$. In the case of (conventional) vector
space, if $y\in Y$, then $Y\oplus y=Y$. However, in the case of 
nonstandard vector space, we just have $Y\oplus y\subseteq Y$ and 
$Y\oplus\omega\subseteq Y\oplus y$, where $\omega =y\ominus y\in\Omega$, 
since, for any $\bar{y}\in Y$, we have $\bar{y}\oplus\omega =
(\bar{y}\ominus y)\oplus y\in Y\oplus y$.
}\end{Rem}

For any $x\in X$, since the distributive law $(\alpha +\beta )x=
\alpha x\oplus\beta x$ does not hold true in general as shown in Example~\ref{nmsex10}, 
it says that $\alpha_{1}x\oplus\cdots\oplus\alpha_{n}x\neq
(\alpha_{1}+\cdots +\alpha_{n})x$ in general. Therefore, we need to 
carefully interpret the concept of linear combination.
Let $X$ be a nonstandard vector space over $\mathbb{F}$ and 
let $\{x_{1},\cdots ,x_{n}\}$ be a finite subset of $X$.
A {\em linear combination} of $\{x_{1},\cdots ,x_{n}\}$ is an expression 
of the form 
\[y_{1}\oplus\cdots\oplus y_{m},\]
where $y_{i}=\alpha_{i}x_{k}$ for some $k\in\{1,\cdots ,n\}$ 
and the coefficients $\alpha_{i}\in\mathbb{F}$ for $i=1,\cdots ,m$. 
We allow $y_{i}=\alpha_{i}x_{k}$ and $y_{j}=\alpha_{j}x_{k}$ for the 
same $x_{k}$. In this case, we see that 
\[y_{i}\oplus y_{j}=\alpha_{i}x_{k}\oplus
\alpha_{j}x_{k}\neq (\alpha_{i}+\alpha_{j})x_{k}\] 
in general. For any nonempty subset $S$ of $X$, the set of all linear combinations of 
finite subsets of $S$ is called the {\em span} of $S$, which is also 
denoted by $\mbox{span}(S)$. Then, it follows that $S\subseteq\mbox{span}(S)$.

\begin{Rem}{\label{ch1r6}}
{\em
Let $X$ be a nonstandard vector space over $\mathbb{R}$.
For $x\in X$, we see that each element in $\mbox{span}(\{x\})$ has the form of
$\alpha_{1}x\oplus\cdots\oplus\alpha_{n}x$ for some finite sequence 
$\{\alpha_{1},\cdots ,\alpha_{n}\}$ in $\mathbb{R}$.
Now suppose that the distributive law for positive and negative
scalar addition hold true in $X$, i.e., 
$(\alpha +\beta )x=\alpha x\oplus\beta x$ for any 
$x\in X$ and $\alpha\beta >0$. For example, 
the distributive law for positive and negative
scalar addition hold true in the set ${\cal I}$ of all closed intervals.
Let $J_{0}=\{j:\alpha_{j}=0\}$,
$J_{+}=\{j:\alpha_{j}>0\}$ and $J_{-}=\{j:\alpha_{j}<0\}$.
We also write $\alpha^{+}=\sum_{j\in J_{+}}\alpha_{j}$ and 
$\alpha^{-}=\sum_{j\in J_{-}}\alpha_{j}$. Then, for $J_{0}\neq\emptyset$,
we see that
\[\alpha_{1}x\oplus\cdots\oplus\alpha_{n}x=\left\{\begin{array}{ll}
\alpha^{+}x\oplus\alpha^{-}x\oplus 0x\oplus\cdots\oplus 0x & 
\mbox{if $J_{+}\neq\emptyset$ and $J_{-}\neq\emptyset$}\\
\alpha^{+}x\oplus 0x\oplus\cdots\oplus 0x & 
\mbox{if $J_{+}\neq\emptyset$ and $J_{-}=\emptyset$}\\
\alpha^{-}x\oplus 0x\oplus\cdots\oplus 0x & 
\mbox{if $J_{+}=\emptyset$ and $J_{-}\neq\emptyset$}\\
0x\oplus\cdots\oplus 0x & 
\mbox{if $J_{+}=\emptyset$ and $J_{-}=\emptyset$}
\end{array}\right .\]
For $J_{0}=\emptyset$, we can have the similar expression without 
$0x\oplus\cdots\oplus 0x$. We need to remark that $0x$ is not necessarily in the 
null set $\Omega$. However, we have the following relations. Suppose that the 
distributive law $(0+0)x=0x+0x$ holds true in general. Then, we have 
$0x=(0+0)x=0x+0x$ by adding $-0x$ on both sides, we have $\omega =0x+\omega$,
where $\omega =0x-0x\in\Omega$.
}\end{Rem}

\begin{Pro}{\label{ch1l6}}
Let $X$ be a nonstandard vector space over $\mathbb{F}$ such that the 
distributive law for vector addition and the associative law for 
scalar multiplication hold true. Then $\mbox{\em span}(S)$ is a subspace of $X$.
\end{Pro}

\section{Decomposition}

In the (conventional) vector space $X$, any element $x\in X$ can be 
decomposed as $x=\hat{x}+(x-\hat{x})$ for some $\hat{x}\in X$. However,
under the nonstandrd vector space $X$, we cannot have the decomposition 
$x=\hat{x}\oplus (x\ominus\hat{x})$ in 
general, since, in fact, $\hat{x}\oplus (x\ominus\hat{x})=x\oplus\omega\neq
x$, where $\omega =\hat{x}\ominus\hat{x}$. 
Now, we propose a very basic notion of decomposition.

\begin{Def}{\label{ch1d556}}
{\em
Let $X$ be a nonstandard vector space over $\mathbb{F}$. 
Given any $x\in X$, we say that $x$ has the {\em null decomposition} 
if $x=\bar{x}\oplus\omega$ for some $\bar{x}\in X$ and $\omega\in\Omega$.
The space $X$ is said to own the null decomposition 
if every element $x\in X$ has the null decomposition. 
}\end{Def}

\begin{Rem}
{\em
Let $X$ be a nonstandard vector space with the null set $\Omega$. Then 
it is easy to see that $X\oplus\Omega\subseteq X$. Now if $X$ owns the 
null decomposition, then we see that $X=X\oplus\Omega$.
}\end{Rem}

The following example shows that the null decomposition
is automatically satisfied in the 
nonstandard vector space of all closed intervals.

\begin{Ex}{\label{ch1e271}}
{\em
Let ${\cal I}$ be the nonstandard vector space of all closed intervals 
with the null set $\Omega =\{[-k,k]:k\geq 0\}$.
For $x=[a,b]\in {\cal I}$, 
we can take $k<(b-a)/2$ and $\bar{x}=[a+k,b-k]$. In this case, 
$x=[a,b]=[a+k,b-k]\oplus [-k,k]=\bar{x}\oplus\omega$. Therefore, ${\cal I}$ owns
the null decomposition.
}\end{Ex}

\begin{Def}
{\em
Let $X$ be a nonstandard vector space over $\mathbb{F}$.
\begin{itemize}
\item Given a fixed $\omega_{0}\in\Omega$, we say that $\Omega$ owns the 
{\em self-decomposition with respect to $\omega_{0}$} 
if every $\omega\in\Omega$ can be represented 
as $\omega =\omega^{\prime}\oplus\omega_{0}$ 
for some $\omega^{\prime}\in\Omega$. 

\item We say that $\Omega$ owns the {\em self-decomposition} if $\Omega$ 
owns the self-decomposition with respect to every $\omega_{0}\in\Omega$. 
\end{itemize}
}\end{Def}

\begin{Ex}{\label{ch1ex388}}
{\em
Let ${\cal I}$ be the nonstandard vector space of all closed intervals 
with the null set $\Omega =\{[-k,k]:k\geq 0\}$. Given 
$\omega =[-k,k]$ for $k\neq 0$, i.e., $k>0$, we can write $k=k_{1}+k_{2}$ 
with $k_{1},k_{2}>0$. Then, we have 
\[\omega =[-k,k]=[-(k_{1}+k_{2}),k_{1}+k_{2}]=
[-k_{1},k_{1}]\oplus [-k_{2},k_{2}],\]
where $\omega_{1}=[-k_{1},k_{1}],\omega_{2}=[-k_{2},k_{2}]\in\Omega$. 
This says that $\Omega\oplus\Omega =\Omega$ under the space ${\cal I}$.
In other words, the null set $\Omega$ owns the self-decomposition.
}\end{Ex}

\section{Nonstandard Metric Spaces}

Now, we are in a position to introduce the concept of 
the so-called nonstandard metric space.

\begin{Def}{\label{ch1d7}}
{\em
Let $X$ be a nonstandard vector space over $\mathbb{F}$
with the null set $\Omega$. 
For the nonnegative real-valued function $d$ defined on $X\times X$,
we consider the following conditions:
\begin{enumerate}
\item [(i)] $d(x,y)=0$ if and only if $x\stackrel{\Omega}{=}y$ for all 
$x,y\in X$;

\item [(i')] $d(x,y)=0$ if and only if $x=y$ for all $x,y\in X$;

\item [(ii)] $d(x,y)=d(y,x)$ for all $x,y\in X$;

\item [(iii)] $d(x,y)\leq d(x,z)+d(z,y)$ for all $x,y,z\in X$;
\end{enumerate}
Different kinds of metric spaces are defined below.
\begin{itemize}
\item A pair $(X,d)$ is called a {\em pseudo-metric space on a nonstandard vector 
space $X$} if and only if $d$ satisfies conditions (ii) and (iii).

\item A pair $(X,d)$ is called a {\em metric space  on a nonstandard vector space $X$} 
if and only if $d$ satisfies conditions (i'), (ii) and (iii).

\item A pair $(X,d)$ is called a {\em nonstandard metric space} if and only if 
$d$ satisfies conditions (i), (ii) and (iii).
\end{itemize}
Now we consider the following conditions:
\begin{enumerate}
\item [(iv)] for any $\omega_{1},\omega_{2}\in\Omega$ and $x,y,z\in X$, we have 
\[d(x\oplus\omega_{1},y\oplus\omega_{2})\geq d(x,y),\quad
d(x\oplus\omega_{1},y)\geq d(x,y)
\mbox{ and }d(x,y\oplus\omega_{2})\geq d(x,y);\]

\item [(iv')] for any $\omega_{1},\omega_{2}\in\Omega$ and $x,y,z\in X$, we have
\[d(x\oplus\omega_{1},y\oplus\omega_{2})=d(x,y),\quad
d(x\oplus\omega_{1},y)=d(x,y)\mbox{ and }d(x,y\oplus\omega_{2})=d(x,y).\]
\end{enumerate}
We say that $d$ satisfies the {\em null inequalities} if and only if condition (iv) 
is satisfied, and that $d$ satisfies the {\em null equalities} if and only if
condition (iv') is satisfied. We also consider the following conditions:
\begin{enumerate}
\item [(v)] for any $x,y,a,b\in X$ and any finite sequences
$\{\alpha_{1},\cdots ,\alpha_{n}\}$ and $\{\beta_{1},\cdots ,\beta_{m}\}$
in $\mathbb{F}$ with $\sum_{i=1}^{n}\alpha_{i}=0$ and 
$\sum_{j=1}^{m}\beta_{j}=0$, we have 
\begin{eqnarray*}
&& d(x\oplus\alpha_{1}a\oplus\cdots\oplus\alpha_{n}a,
y\oplus\beta_{1}b\oplus\cdots\oplus\beta_{m}b)\geq d(x,y),\\
&& d(x\oplus\alpha_{1}a\oplus\cdots\oplus\alpha_{n}a,y)\geq d(x,y)\mbox{ and }\\
&& d(x,y\oplus\beta_{1}b\oplus\cdots\oplus\beta_{m}b)\geq d(x,y);
\end{eqnarray*}

\item [(v')] for any $x,y,a,b\in X$ and any finite sequences
$\{\alpha_{1},\cdots ,\alpha_{n}\}$ and $\{\beta_{1},\cdots ,\beta_{m}\}$
in $\mathbb{F}$ with $\sum_{i=1}^{n}\alpha_{i}=0$ and 
$\sum_{j=1}^{m}\beta_{j}=0$, we have 
\begin{eqnarray*}
&& d(x\oplus\alpha_{1}a\oplus\cdots\oplus\alpha_{n}a,
y\oplus\beta_{1}b\oplus\cdots\oplus\beta_{m}b)=d(x,y),\\
&& d(x\oplus\alpha_{1}a\oplus\cdots\oplus\alpha_{n}a,y)=d(x,y)\mbox{ and }\\
&& d(x,y\oplus\beta_{1}b\oplus\cdots\oplus\beta_{m}b)=d(x,y).
\end{eqnarray*}
\end{enumerate}
We say that the metric (or pseudo-metric) $d$ satisfies the 
{\em zero-sum inequalities} if and only if condition (v) is satisfied,
and that $d$ satisfies the {\em zero-sum equalities} if and only if 
condition (v') is satisfied.
}\end{Def}

\begin{Rem}{\label{ch1*r15}}
{\em
For any $\omega_{1},\omega_{2}\in\Omega$, since $\omega_{1}=a\ominus a$
and $\omega_{2}=b\ominus b$ for some $a,b\in X$, it is not hard to 
see that condition (v) implies condition (iv), and 
condition (v') implies condition (iv'). 
}\end{Rem}

\begin{Ex}{\label{ch1ex200}}
{\em
We consider the nonstandard vector space ${\cal I}$ of all
closed intervals in $\mathbb{R}$. Let us define a nonnegative 
real-valued function $d:{\cal I}\times {\cal I}
\rightarrow\mathbb{R}_{+}$ by 
\begin{equation}{\label{ch1eq203}}
d([a,b],[c,d])=|(a+b)-(c+d)|.
\end{equation}
Then we are going to claim that $({\cal I},d)$ is a nonstandard 
metric space such that the metric $d$ satisfies the zero-sum equalities.
Remark~\ref{ch1*r15} says that the metric also satisfies the null equalities.
\begin{itemize}
\item We consider the closed intervals $[a,b]$ and $[c,d]$. 
Then we see that 
$a-d\leq b-c$. Therefore, if $b-c<0$, then $d([a,b],[c,d])=|a+b-c-d|\neq 0$.
Suppose that $b-c\geq 0$ and $d([a,b],[c,d])=|a+b-c-d|=0$. Then we have 
$a+b=c+d$, i.e., $a+c-d=2c-b$. It is easy to see that $a+c-d\leq b+d-c$
and $2c-b\leq b+d-c$ by using the facts that $a\leq b$, $c\leq d$ and 
$b\geq c$ (in fact, this can be understood from (\ref{ch1eq63}) below). 
Therefore we can form the two idetntical closed intervals 
$[a+c-d,b+d-c]=[2c-b,b+d-c]$. Now the closed intervals 
$[a+c-d,b+d-c]$ and $[2c-b,b+d-c]$ can be written as 
\begin{equation}{\label{ch1eq63}}
[a+c-d,b+d-c]=[a,b]\oplus [c-d,d-c]\mbox{ and }[2c-b,b+d-c]=
[c,d]\oplus [c-b,b-c].
\end{equation}
Let $\omega_{1}=[c-d,d-c]=(d-c)[-1,1]\in\Omega$ and 
$\omega_{2}=[c-b,b-c]=(b-c)[-1,1]\in\Omega$. Therefore, from (\ref{ch1eq63}),
we have $[a,b]\oplus\omega_{1}=[c,d]\oplus\omega_{2}$,
which shows $[a,b]\stackrel{\Omega}{=} [c,d]$, since $\omega_{1},\omega_{2}\in\Omega$.
Conversely, suppose that $[a,b]\stackrel{\Omega}{=} [c,d]$. Then 
$[a,b]\oplus\omega_{1}=[c,d]\oplus\omega_{2}$,
where $\omega_{1}=[-k_{1},k_{1}],\omega_{2}=[-k_{2},k_{2}]\in\Omega$.
Therefore, we have $[a-k_{1},b+k_{1}]=[c-k_{2},d+k_{1}]$, i.e., 
$a-k_{1}=c-k_{2}$ and $b+k_{1}=d+k_{2}$. Then we obtain
\[d([a,b],[c,d])=|(a-c)+(b-d)|=|(k_{1}-k_{2})+(k_{2}-k_{1})|=0.\]

\item We have
\[d([a,b],[c,d])=|a+b-c-d|=|c+d-a-b|=d([c,d],[a,b]).\]

\item We have 
\begin{align*}
d([a,b],[c,d]) & =|a+b-c-d|=|(a+b-e-f)+(e+f-c-d)|\\
& \leq |a+b-e-f|+|e+f-c-d|\\
& =d([a,b],[e,f])+d([e,f],[c,d]).
\end{align*}

\item For any finite sequence $\{\alpha_{1},\cdots ,\alpha_{n}\}$ in 
$\mathbb{R}$ with $\sum_{i=1}^{n}\alpha_{i}=0$, we have
\[\alpha_{1}[c_{1},d_{1}]\oplus\cdots\oplus\alpha_{n}[c_{1},d_{1}]
=[e_{1},f_{1}],\]
where 
\[e_{1}=\sum_{\alpha_{i}\geq 0}\alpha_{i}c_{1}+\sum_{\alpha_{i}<0}\alpha_{i}d_{1}
\mbox{ and }f_{1}=\sum_{\alpha_{i}\geq 0}\alpha_{i}d_{1}+
\sum_{\alpha_{i}<0}\alpha_{i}c_{1}.\]
Therefore
\[e_{1}+f_{1}=\sum_{\alpha_{i}\geq 0}\alpha_{i}(c_{1}+d_{1})+\sum_{\alpha_{i}<0}
\alpha_{i}(c_{1}+d_{1})=(c_{1}+d_{1})\sum_{i=1}^{n}\alpha_{i}=0.\]
Similarly, for any finite sequence $\{\beta_{1},\cdots ,\beta_{m}\}$ in 
$\mathbb{R}$ with $\sum_{j=1}^{m}\beta_{j}=0$, we have
$\beta_{1}[c_{2},d_{2}]\oplus\cdots\oplus\beta_{m}[c_{2},d_{2}]=[e_{2},f_{2}]$ 
implies $e_{2}+f_{2}=0$. Now, we have
\begin{align*}
& d([a_{1},b_{1}]\oplus\alpha_{1}[c_{1},d_{1}]\oplus\cdots\oplus\alpha_{n}[c_{1},d_{1}],
[a_{2},b_{2}]\oplus\beta_{1}[c_{2},d_{2}]\oplus\cdots\oplus\beta_{m}[c_{2},d_{2}])\\
& \quad =d([a_{1},b_{1}]\oplus [e_{1},f_{1}],[a_{2},b_{2}]\oplus [e_{2},f_{2}])\\
& \quad =d([a_{1}+e_{1},b_{1}+f_{1}],[a_{2}+e_{2},b_{2}+f_{2}])\\
& \quad =|(a_{1}+e_{1}+b_{1}+f_{1})-(a_{2}+e_{2}+b_{2}+f_{2})|\\
& \quad =|(a_{1}+b_{1})-(a_{2}+b_{2})|\\
& \quad =d([a_{1},b_{1}],[a_{2},b_{2}]).
\end{align*}
\end{itemize}
This shows that $({\cal I},d)$ 
is also a nonstandard metric space such that the metric 
$d$ satisfies the zero-sum and null equalities.
We also remark that $({\cal I},d)$ cannot be a metric space,
since condition (i') in Definition~\ref{ch1d7} cannot hold true.
}\end{Ex}

\begin{Ex}{\label{ch1ex143}}
{\em
Example~\ref{ch1ex200} shows that $({\cal I},d_{\cal I})$ 
is a nonstandard metric 
space such that the metric $d$ satisfies the zero-sum equalities define in 
(\ref{ch1eq203}). Let $(X,d)$ be another nonstandard 
metric space such that the metric $d$ satisfies the zero-sum equalities. 
We consider the interval-valued function $F$ defined on 
$X$. In other words, for each $x\in X$, $F(x)=[f^{L}(x),f^{U}(x)]$ is a
closed interval in $\mathbb{R}$. Then we say that $F$ is continuous at 
$x_{0}$ if, for every $\epsilon >0$, there exists $\delta >0$ such 
that $d(x,x_{0})<\delta$ implies $d_{\cal I}(F(x),F(x_{0}))<\epsilon$.
We see that 
\[d_{\cal I}(F(x),F(x_{0}))=\left |f^{L}(x)+f^{U}(x)-f^{L}(x_{0})-f^{U}(x_{0})
\right |.\]
We also say that $F$ is continuous on $X$ if $F$ is continuous at 
each $x_{0}\in X$. We also denote by ${\cal IC}(X)$ the set of all continuous
functions $F:(X,d)\rightarrow (I,d_{\cal I})$ on $X$. Now the vector addition 
and scalcar multiplications in ${\cal IC}(X)$ are defined by
\[(F\oplus G)(x)=F(x)\oplus G(x)\mbox{ and }(\alpha F)(x)=\alpha F(x)\]
for any $F,G\in {\cal IC}(X)$ and $\alpha\in \mathbb{R}$. Then we are going to show 
that $F\oplus G$ and $\alpha F$ are also in ${\cal IC}(X)$. Since 
\begin{align*}
& d_{\cal I}((F\oplus G)(x),(F\oplus G)(x_{0}))\\
& \quad =\left |f^{L}(x)+f^{U}(x)+g^{L}(x)+
g^{U}(x)-f^{L}(x_{0})-f^{U}(x_{0})-g^{L}(x_{0})-g^{U}(x_{0})\right |\\
& \quad\leq\left |f^{L}(x)+f^{U}(x)-f^{L}(x_{0})-f^{U}(x_{0})\right |+
\left |g^{L}(x)+g^{U}(x)-g^{L}(x_{0})-g^{U}(x_{0})\right |
\end{align*}
and
\[d_{\cal I}(\alpha F(x),\alpha F(x_{0}))
=|\alpha |\cdot |f^{L}(x)+f^{U}(x)-f^{L}(x_{0})-f^{U}(x_{0})|,\]
we see that $F\oplus G$ and $\alpha F$ are continuous on $X$, since 
$F$ and $G$ are continuous on $X$. The commutative and
associative laws for vector addition hold true obviously. 
Therefore, we conclude that ${\cal IC}(X)$ is a nonstandard vector space.
It is not hard to see that the null set us
\[\Omega_{\cal IC}=\{[-k(x),k(x)]:k(x)\geq 0\mbox{ for all }x\in X\}.\]
Now we want to introduce a metric $d_{\cal IC}$ to make 
$({\cal IC}(X),d_{\cal IC})$ as a
nonstandard metric space such that the metric $d_{\cal IC}$
satisfies the zero-sum equalities. For $F,G\in {\cal IC}(X)$, we define
\[d_{\cal IC}(F,G)=\sup_{x\in X}d_{\cal I}(F(x),G(x)).\]
We need to check four conditions.
\begin{itemize}
\item We have that 
\[0=d_{\cal IC}(F,G)=\sup_{x\in X}d_{\cal I}(F(x),G(x))\]
implies $d_{\cal I}(F(x),G(x))=0$ for all $x\in X$, i.e., $F(x)\stackrel{\Omega}{=} G(x)$
for all $x\in X$, which also says that $F\stackrel{\Omega}{=} G$ in the sense of 
$\Omega_{\cal IC}$

\item Since $d_{\cal I}$ is symmetric, it is easy to see that $d_{\cal IC}$ is symmetric. 

\item We have 
\begin{align*}
d_{\cal IC}(F,G) & =\sup_{x\in X}d_{\cal I}(F(x),G(x))\leq
\sup_{x\in X}\left [d_{\cal I}(F(x),E(x))+d_{\cal I}(E(x),G(x))\right ]\\
& \leq\sup_{x\in X}d_{\cal I}(F(x),E(x))+\sup_{x\in X}d_{\cal I}(E(x),G(x))
=d_{\cal IC}(F,E)+d_{\cal IC}(E,G).
\end{align*}

\item For any finite sequences $\{\alpha_{1},\cdots ,\alpha_{n}\}$ and
$\{\beta_{1},\cdots ,\beta_{m}\}$ in $\mathbb{R}$ with 
$\sum_{i=1}^{n}\alpha_{i}=0$ and $\sum_{j=1}^{m}\beta_{j}=0$ and 
$H_{1},H_{2}\in {\cal IC}(X)$, we see that
$\alpha_{1}H_{1}(x)\oplus\cdots\oplus\alpha_{n}H_{1}(x)$
and $\beta_{1}H_{2}(x)\oplus\cdots\oplus\beta_{m}H_{2}(x)$
are also in ${\cal IC}(X)$. Therefore, we have 
\begin{align*}
& d_{\cal IC}(F\oplus\alpha_{1}H_{1}\oplus\cdots\oplus\alpha_{n}H_{1},
G\oplus\beta_{1}H_{2}\oplus\cdots\oplus\beta_{m}H_{2})\\
& \quad =\sup_{x\in X}d_{\cal I}(F(x)\oplus\alpha_{1}H_{1}(x)\oplus\cdots\oplus\alpha_{n}H_{1}(x),
G(x)\oplus\beta_{1}H_{2}(x)\oplus\cdots\oplus\beta_{m}H_{2}(x))\\
& \quad =\sup_{x\in X}d_{\cal I}(F(x),G(x))=d_{\cal IC}(F,G).
\end{align*}
\end{itemize}
Therefore, we conclude that $({\cal IC}(X),d_{\cal IC})$ is indeed a nonstandard 
metric space such that the metric $d_{\cal IC}$ satisfies the zero-sum equalities.
}\end{Ex}

\begin{Def}
{\em
Let $(X,d)$ be a pseudo-metric space on a nonstandard vector space $X$. 
We say that the pseudo-metric $d$ is {\em translation-invariant} if and only if 
$d(x\oplus z,y\oplus z)=d(x,y)$. We say that $d$ is {\em absolutely homogeneous} 
if and only if $d(\alpha x,\alpha y)=|\alpha |d(x,y)$.
}\end{Def}

\section{Open and Closed Balls}

Let $(X,d)$ be a pseudo-metric space on a nonstandard vector space $X$.
Given a point $x_{0}\in X$ and a positive number $\epsilon>0$, the 
{\em open ball} about $x_{0}$ is defined by
\[B(x_{0};\epsilon )=\{x\in X:d(x,x_{0})<\epsilon\},\]
the {\em closed ball} about $x_{0}$ is defined by 
\[\bar{B}(x_{0};\epsilon )=\{x\in X:d(x,x_{0})\leq\epsilon\}\]
and the {\em sphere} about $x_{0}$ is defined by 
\[S(x_{0};\epsilon )=\{x\in X:d(x,x_{0})=\epsilon\}.\]
In all three cases, $x_{0}$ is called the {\em center} and $\epsilon$ the 
{\em radius}. 

\begin{Pro}{\label{ch2p91}}
Let $X$ be a nonstandard vector space over $\mathbb{F}$, and let 
$(X,d)$ be a pseudo-metric space. The following statements hold true.
\begin{enumerate}
\item [{\em (i)}] If $d$ satisfies the null inequalities, then 
$x\oplus\omega\in B(x_{0};\epsilon )$ implies $x\in B(x_{0};\epsilon )$
for any $\omega\in\Omega$.

\item [{\em (ii)}] If $d$ satisfies the null equalities, then
$x\oplus\omega\in B(x_{0};\epsilon )$ if and only if $x\in B(x_{0};\epsilon )$
for any $\omega\in\Omega$. Moreover, we have the following inclusion 
$B(x_{0};\epsilon )\oplus\Omega\subseteq B(x_{0};\epsilon )$.
\end{enumerate}
\end{Pro}
\begin{Proof}
To prove part (i), suppose that $x\oplus\omega\in B(x_{0};\epsilon )$. 
According to condition (iv) of Definition~\ref{ch1d7}, we have 
\[d(x,x_{0})\leq d(x\oplus\omega ,x_{0})<\epsilon ,\]
which shows that $x\in B(x_{0};\epsilon )$. To prove part (ii), 
for $x\in B(x_{0};\epsilon )$ and $\omega\in\Omega$, 
according to condition (iv') of Definition~\ref{ch1d7}, it follows that 
\[d(x_{0},x\oplus\omega )=d(x_{0},x)<\epsilon ,\]
which says that $x\oplus\omega\in B(x_{0};\epsilon )$. This shows the
desired inclusion, and the proof is complete.
\end{Proof}

\begin{Pro}{\label{ch2p1}}
Let $(X,d)$ be a pseudo-metric space on a nonstandard vector space $X$. 
The following statements hold true.
\begin{enumerate}
\item [{\em (i)}] If the pseudo-metric $d$ satisfies the null inequalities, then 
$B(x_{0}\oplus\omega ;\epsilon )\subseteq B(x_{0};\epsilon )$ 
for any $\omega\in\Omega$.

\item [{\em (ii)}] If the pseudo-metric $d$ satisfies the null equalities, then 
$B(x_{0}\oplus\omega ;\epsilon )=B(x_{0};\epsilon )$ for any $\omega\in\Omega$.
\end{enumerate}
\end{Pro}
\begin{Proof}
If the pseudo-metric $d$ satisfies the null inequalities, then the inclusion 
$B(x_{0}\oplus\omega ;\epsilon )\subseteq B(x_{0};\epsilon )$ follows from 
the ineqaulity $\epsilon >d(x,x_{0}\oplus\omega)\geq d(x,x_{0})$ immediately. 
Suppose that the pseudo-metric $d$ satisfies the null equalities. Then the inclusion
$B(x_{0};\epsilon )\subseteq B(x_{0}\oplus\omega ;\epsilon )$ follows from the 
equality $\epsilon >d(x,x_{0})=d(x,x_{0}\oplus\omega)$ immediately.
This completes the proof.
\end{Proof}

\begin{Pro}{\label{ch1p273}}
Let $X$ be a nonstandard vector space over $\mathbb{F}$, and let 
$(X,d)$ be a pseudo-metric space. Suppose that $X$ owns the null decomposition.
The following statements hold true.
\begin{enumerate}
\item [{\em (i)}] If $d$ satisfies the null inequalities,
then $B(x_{0};\epsilon )\subseteq B(x_{0};\epsilon )\oplus\Omega$.

\item [{\em (ii)}] If $d$ satisfies the null equalities, 
then $B(x_{0};\epsilon )\oplus\Omega =B(x_{0};\epsilon )$.
\end{enumerate}
\end{Pro}
\begin{Proof}
For any $x\in B(x_{0};\epsilon )$, since $x=\bar{x}\oplus\omega$ for some 
$\bar{x}\in X$ and $\omega\in\Omega$ by the null decomposition, we have 
\[d(\bar{x},x_{0})\leq d(\bar{x}\oplus\omega,x_{0})=d(x,x_{0})<\epsilon .\]
This says that $\bar{x}\in B(x_{0};\epsilon )$, i.e., $x=\bar{x}\oplus
\omega\in B(x_{0};\epsilon )\oplus\Omega$. This proves part (i). 
Part (ii) follows from part (ii) of Proposition~\ref{ch2p91} and part (i).
This completes the proof.
\end{Proof}

In the (conventional) pseudo-metric space $(X,d)$, where $X$ is taken as a 
(conventional) vector space. If $d$ is translation-invariant, 
then we have the following equality 
\begin{equation}{\label{nnaeq3}}
B(x;\epsilon )\oplus\hat{x}=B(x\oplus\hat{x};\epsilon ).
\end{equation}
However, if $X$ is taken as the nonstandard vector space, then 
the intuitive observation in (\ref{nnaeq3}) will not hold true.
The following proposition presents the exact relationship, and will be used 
for studying the topology induced by the nonstandard metric space.

\begin{Pro}{\label{ch1r289}}
Let $X$ be a nonstandard vector space over $\mathbb{F}$, and let 
$(X,d)$ be a pseudo-metric space. The following statements hold true.
\begin{enumerate}
\item [{\em (i)}] If $d$ is translation-invariant, then 
\begin{equation}{\label{nnaeq5}}
B(x;\epsilon )\oplus\hat{x}\subseteq B(x\oplus\hat{x};\epsilon )
\mbox{ for any }\omega\in\Omega .
\end{equation} 
and
\begin{equation}{\label{nnaeq6}}
B(x;\epsilon )\oplus\omega_{x}\subseteq x\oplus B(\omega_{x};\epsilon ),
\mbox{ where }\omega_{x}=x\ominus x\in\Omega .
\end{equation} 

\item [{\em (ii)}] If $d$ is translation-invariant and satisfies the null inequalities, 
then 
\begin{equation}{\label{nnaeq7}}
B(x;\epsilon )\oplus\omega\subseteq B(x;\epsilon )\mbox{ and }
B(\omega ;\epsilon )\oplus\hat{x}\subseteq B(\hat{x};\epsilon )
\mbox{ for any }\omega\in\Omega .
\end{equation} 
and
\begin{equation}{\label{nnaeq8}}
B(x\oplus\hat{x};\epsilon )\oplus\omega_{\hat{x}}\subseteq 
B(x;\epsilon )\oplus\hat{x},
\mbox{ where }\omega_{\hat{x}}=\hat{x}\ominus\hat{x}.
\end{equation} 

\item [{\em (iii)}] Suppose that $d$ is translation-invariant
and satisfies the null inequalities.
We also assume that $X$ owns the null decomposition and $\Omega$ owns 
the self-decomposition, then we have the following equality
\begin{equation}{\label{nnaeq9}}
B(x\oplus\hat{x};\epsilon )=B(x;\epsilon )\oplus\hat{x}
\mbox{ for any }\omega\in\Omega .
\end{equation}

\item [{\em (iv)}] Suppose that $d$ is translation-invariant
and satisfies the null equalities.
We also assume that $X$ owns the null decomposition and $\Omega$ owns 
the self-decomposition, then we have the following equalities 
\[B(x;\epsilon )\oplus\omega =B(x;\epsilon )\mbox{ and }
B(\omega ;\epsilon )\oplus\hat{x}=B(\hat{x};\epsilon )\]
for any $\omega\in\Omega$.
\end{enumerate}
\end{Pro}
\begin{Proof}
To prove part (i), for $y\in B(x;\epsilon )\oplus\hat{x}$,
we have $y=\hat{y}\oplus\hat{x}$ with $d(\hat{y},x)<\epsilon$. 
Then we can obtain
\[d(y,x\oplus\hat{x})=d(\hat{y}\oplus\hat{x},x\oplus\hat{x})=
d(\hat{y},x)<\epsilon ,\]
which says that $y\in B(x\oplus\hat{x};\epsilon )$. Therefore, we obtain 
(\ref{nnaeq5}). To prove (\ref{nnaeq6}),
for $\hat{x}\in B(x;\epsilon )$ and $\omega_{x}=x\ominus x$, 
we have $\hat{x}\oplus\omega_{x}=x\oplus (\hat{x}\ominus x)$. Then we obtain
\[d(\hat{x}\ominus x,\omega_{x})=d(\hat{x}\ominus x,x\ominus x)=
d(\hat{x}\oplus (-x),x\oplus (-x))=d(\hat{x},x)<\epsilon ,\]
which says that $\hat{x}\ominus x\in B(\omega_{x};\epsilon )$. This shows that 
$\hat{x}\oplus\omega_{x}=x\oplus (\hat{x}\ominus x)\in x\oplus 
B(\omega_{x};\epsilon )$.

To prove part (ii), we take $\hat{x}=\omega\in\Omega$ in (\ref{nnaeq5}).
By part (i) of Proposition~\ref{ch2p1}, we have 
\[B(x;\epsilon )\oplus\omega\subseteq
B(x\oplus\omega ;\epsilon )\subseteq B(x;\epsilon )\] 
for any $\omega\in\Omega$. Similarly, if we take $x=\omega$, then we have 
\[B(\omega ;\epsilon )\oplus\hat{x}\subseteq 
B(\omega\oplus\hat{x};\epsilon )=B(\hat{x};\epsilon ).\]
Therefore, we obtain (\ref{nnaeq7}). To prove (\ref{nnaeq8}),
for $y\in B(x\oplus\hat{x};\epsilon )$, we have 
$d(y,x\oplus\hat{x})<\epsilon$. Since 
\[\epsilon >d(y,x\oplus\hat{x})=d(y\ominus\hat{x},x\oplus\hat{x}\ominus
\hat{x})=d(y\ominus\hat{x},x\oplus\omega_{\hat{x}})\geq d(y\ominus\hat{x},x),\]
we see that $y\ominus\hat{x}\in B(x;\epsilon )$, where $\omega_{\hat{x}}=\hat{x}
\ominus\hat{x}\in\Omega$. Since $y\oplus\omega_{\hat{x}}=\hat{x}\oplus (y\ominus
\hat{x})$, it says that $y\oplus\omega_{\hat{x}}\in B(x;\epsilon )\oplus\hat{x}$,
and proves (\ref{nnaeq8}).

To prove part (iii), 
for $y\in B(x\oplus\hat{x};\epsilon )$, we have $d(y,x\oplus\hat{x})<\epsilon$.
Since $X$ owns the null decomposition, we have $y=\hat{y}\oplus\hat{\omega}_{0}$ 
for some $\hat{y}\in X$ and $\hat{\omega}_{0}\in\Omega$. Since $\Omega$ owns 
the self-decompsition, we also have $\hat{\omega}_{0}=\omega_{\hat{x}}
\oplus\hat{\omega}_{1}$ 
for some $\hat{\omega}_{1}\in\Omega$. Then we have $y=\hat{y}\oplus\omega_{\hat{x}}
\oplus\hat{\omega}_{1}$ and 
\[d(\hat{y}\oplus\hat{\omega}_{1},x\oplus\hat{x})\leq
d(\hat{y}\oplus\hat{\omega}_{1}\oplus\omega_{\hat{x}},x\oplus\hat{x})=
d(y,x\oplus\hat{x})<\epsilon .\]
Since 
\[\epsilon >d(\hat{y}\oplus\hat{\omega}_{1},x\oplus\hat{x})
=d(\hat{y}\oplus\hat{\omega}_{1}\ominus\hat{x},x\oplus\hat{x}\ominus\hat{x})
=d(\hat{y}\oplus\hat{\omega}_{1}\ominus\hat{x},x\oplus\omega_{\hat{x}})
\geq d(\hat{y}\oplus\hat{\omega}_{1}\ominus\hat{x},x),\]
it follows that $\hat{y}\oplus\hat{\omega}_{1}\ominus\hat{x}\in B(x;\epsilon )$. 
By adding $\hat{x}$ on both sides, we have 
\[y=\hat{y}\oplus\hat{\omega}_{1}\oplus\omega_{\hat{x}}
\in B(x;\epsilon )\oplus\hat{x};\]
that is, we have the inclusion 
\[B(x\oplus\hat{x};\epsilon )\subseteq B(x;\epsilon )\oplus\hat{x}.\] 
From (\ref{nnaeq5}), we obtain
\[B(x\oplus\hat{x};\epsilon )=B(x;\epsilon )\oplus\hat{x}.\]

To prove part (iv), we take $\hat{x}=\omega\in\Omega$ in (\ref{nnaeq9}).
By part (ii) of Proposition~\ref{ch2p1}, we have 
\[B(x;\epsilon )\oplus\omega =B(x\oplus\omega ;\epsilon )=B(x;\epsilon )\] 
for any $\omega\in\Omega$. Similarly, if we take $x=\omega$, then we have 
\[B(\omega ;\epsilon )\oplus\hat{x}=
B(\omega\oplus\hat{x};\epsilon )=B(\hat{x};\epsilon ).\]
This completes the proof.
\end{Proof}

\begin{Pro}{\label{ch1l437}}
Let $(X,d)$ be a pseudo-metric space on a nonstandard vector space $X$.
Suppose that the pseudo-metric $d$ satisfies the 
null equalities and $\Omega$ owns the self-decomposition with respect to 
$\omega_{0}$. Then $B(a;\epsilon )\oplus\omega_{0}
=B(a;\epsilon )\oplus\Omega$ for any open ball $B(a;\epsilon )$. 
Of course, if $\Omega$ owns the 
self-decomposition. Then $B(a;\epsilon )\oplus\omega =
B(a;\epsilon )\oplus\Omega$ for any $\omega\in\Omega$ 
and open ball $B(a;\epsilon )$.
\end{Pro}
\begin{Proof}
We see that $B(a;\epsilon )\oplus\omega_{0}\subseteq B(a;\epsilon )
\oplus\Omega$ for any $\omega\in\Omega$. On the other hand,
for any $y\in B(a;\epsilon )\oplus\Omega$, we have $y=\hat{y}
\oplus\hat{\omega}$ for some $\hat{y}\in B(a;\epsilon )$ and 
$\hat{\omega}\in\Omega$. Since $\Omega$ owns the  
self-decomposition with respect to $\omega_{0}$, 
we have $\hat{\omega}=\omega^{\prime}\oplus\omega_{0}$ for some 
$\omega^{\prime}\in\Omega$. Therefore, we have 
$y=\hat{y}\oplus\hat{\omega}=\hat{y}\oplus
\omega^{\prime}\oplus\omega_{0}$.
Then we obtain $d(\hat{y}\oplus\omega^{\prime},a)=d(\hat{y},a)
<\epsilon$, which says that $\hat{y}\oplus\omega^{\prime}\in 
B(a;\epsilon )$. Therefore, we obtain that 
$y\in B(a;\epsilon )\oplus\omega_{0}$. This shows that 
$B(a;\epsilon )\oplus\Omega =B(a;\epsilon )\oplus\omega_{0}$. 
We complete the proof.
\end{Proof}

\begin{Pro}{\label{ch1p554}}
Let $(X,d)$ be a nonstandard metric space.
Suppose that the pseudo-metric $d$ satisfies the null equalities and 
$\Omega$ owns the self-decomposition. Then, given any $w\in\Omega$,
we have $d(\omega^{\prime},\omega )=0$ for any $\omega^{\prime}\in\Omega$. 
In other words, we have $\Omega\subseteq B(\omega ;\epsilon )$
for any $\omega\in\Omega$ and $\epsilon >0$.
\end{Pro}
\begin{Proof}
Given any $w\in\Omega$,, since $\Omega$ owns the self-decomposition, 
we have $\omega^{\prime}=
\omega\oplus\hat{\omega}$ for some $\hat{\omega}\in\Omega$. 
Since $d$ satisfies the null equalities, we have 
$d(\omega^{\prime},\omega )=d(\omega\oplus\hat{\omega},\omega )
=d(\omega ,\omega )=0$, which shows that $\omega^{\prime}\in 
B(\omega ;\epsilon )$. This completes the proof.
\end{Proof}

\begin{Pro}{\label{ch1p555}}
Let $(X,d)$ be a pseudo-metric space on a nonstandard vector space $X$.
Suppose that the following conditions are satisfied.
\begin{itemize}
\item The associative law for scalar multiplication holds true. 

\item The pseudo-metric $d$ is absolutely homogeneous and satisfies
the null inequalities.

\item The null set $\Omega$ is closed under the scalar multiplication 
and owns the self-decomposition.
\end{itemize}
Then $\alpha B(\omega ;\epsilon )=B(\omega ;|\alpha |\epsilon )$ 
for any $\omega\in\Omega$ and $\alpha\neq 0$.
\end{Pro}
\begin{Proof}
For $x\in B(\omega ;\epsilon )$, since $\alpha\omega\in\Omega$, we have 
\[d(\alpha x,\omega )=d(\alpha x,\omega\oplus\alpha\omega )
=d(\alpha x,\alpha\omega )=|\alpha |d(x,\omega )<|\alpha |\epsilon ,\] 
i.e., $\alpha x\in B(\omega ;|\alpha |\epsilon )$. 
This shows the inclusion $\alpha B(\omega ;\epsilon )\subseteq 
B(\omega ;|\alpha |\epsilon )$. On the other hand, for $x\in B(\omega ;
|\alpha |\epsilon )$, we have $d(x,\omega )<|\alpha |\epsilon$, i.e., 
$1/|\alpha |\cdot d(x,\omega )<\epsilon$. 
Since $d$ is absolutely homogeneous, we have 
$d(x/\alpha ,\omega /\alpha )<\epsilon$. 
Since $\Omega$ owns the self-decomposition, we have 
$\omega /\alpha =\omega\oplus\hat{\omega}$ for some $\hat{\omega}\in\Omega$. 
Therefore, we obtain $d(x/\alpha ,\omega )=d(x/\alpha ,\omega\oplus\hat{\omega})=
d(\alpha /x,\omega /\alpha )<\epsilon$, which says that $x/\alpha\in 
B(\omega ;\epsilon )$. Since $1x=x$ and $(\alpha\beta )x=\alpha (\beta x)$, 
we conclude that $x\in\alpha B(\omega ;\epsilon )$. This completes the proof.
\end{Proof}

\section{Nonstandardly Open Sets}

Let $(X,d)$ be a pseudo-metric space on a nonstandard vector space $X$.
We are going to consider the concept of openness of subsets of $X$.

\begin{Def}{\label{nfad2}}
{\em 
Let $A$ be a subset of a pseudo-metric space $(X,d)$
on a nonstandard vector space $X$.
\begin{itemize}
\item A point $x_{0}\in A$ is said to be a {\em nonstandard interior point} of $A$ if 
there exists $\epsilon >0$ such that $B(x_{0};\epsilon )\subseteq A$. 
The collection of all interior points of $A$ is 
called the {\em nonstandard interior} of $A$ and is denoted by 
$\mbox{int}(A)$.

\item A point $x_{0}\in A$ is said to be a {\em nonstandard type-I interior 
point} of $A$ if there exists $\epsilon >0$ such that 
$B(x_{0};\epsilon )\oplus\Omega\subseteq A$. The collection of all 
nonstandard type-I interior points of $A$ is called the 
{\em nonstandard type-I interior} of $A$ and is 
denoted by $\mbox{int}^{\tiny \mbox{(I)}}(A)$.

\item A point $x_{0}\in A$ is said to be a {\em nonstandard type-II interior 
point} of $A$ if there exists $\epsilon >0$ such that 
$B(x_{0};\epsilon )\subseteq A\oplus\Omega$. The collection of all 
nonstandard type-II interior points of $A$ is called the 
{\em nonstandard type-II interior} of $A$ and is 
denoted by $\mbox{int}^{\tiny \mbox{(II)}}(A)$.

\item A point $x_{0}\in A$ is said to be a {\em nonstandard type-III interior 
point} of $A$ if there exists $\epsilon >0$ such that 
$B(x_{0};\epsilon )\oplus\Omega\subseteq A\oplus\Omega$. The collection of all 
nonstandard type-III interior points of $A$ is called the 
{\em nonstandard type-III interior} of $A$ and is 
denoted by $\mbox{int}^{\tiny \mbox{(III)}}(A)$.
\end{itemize}
}\end{Def}

We see that the concept of interior point in Definition~\ref{nfad2}
is the same as the conventional definition. 
We also remark that if $X$ happens to be a (conventional) vector space, then 
$\Omega =\{\theta\}=$ (a zero element of $X$). In this case, the four 
concepts of (nonstandard) interior points coincide with the conventional definition.

\begin{Rem}{\label{ch1r239}}
{\em
Let $(X,d)$ be a pseudo-metric space on a nonstandard vector space $X$.
We have the following observations
\begin{itemize}
\item Since $B(x_{0};\epsilon )\subseteq A$ implies $B(x_{0};\epsilon )
\oplus\Omega\subseteq A\oplus\Omega$, a nonstandard 
interior point is also a nonstandard type-III interior point.
In other words, we have 
$\mbox{int}(A)\subseteq\mbox{int}^{\tiny \mbox{(III)}}(A)$.

\item If we assume that $\Omega\oplus\Omega =\Omega$, then a nonstandard 
type-I interior point is also a nonstandard type-III interior point, since 
$B(x_{0};\epsilon )\oplus\Omega\subseteq A$ implies 
$B(x_{0};\epsilon )\oplus\Omega\oplus\Omega\subseteq A\oplus\Omega$. 
In other words, we have 
$\mbox{int}^{\tiny \mbox{(I)}}(A)\subseteq\mbox{int}^{\tiny \mbox{(III)}}(A)$.

\item Suppose that the pseudo-metric $d$ satisfies the null equalities. 
From Proposition~\ref{ch2p91}, 
we see that $B(x;\epsilon )\oplus\Omega\subseteq B(x;\epsilon )$.
It says that if $x$ is a nonstandard interior point, 
then it is also a nonstandard type-I 
interior point, and if $x$ is a nonstandard type-II interior point, then it is also a 
nonstandard type-III interior point. In other words, we have 
$\mbox{int}(A)\subseteq\mbox{int}^{\tiny \mbox{(I)}}(A)$ and 
$\mbox{int}^{\tiny \mbox{(II)}}(A)\subseteq\mbox{int}^{\tiny \mbox{(III)}}(A)$.

\item Suppose that the pseudo-metric $d$ satisfies the 
null inequalities and $X$ owns the null decomposition. 
From Proposition~\ref{ch1p273}, 
we see that $B(x;\epsilon )\subseteq B(x;\epsilon )\oplus\Omega$.
It says that if $x$ is a nonstandard type-I interior point, then it is also a
nonstandard interior point, and if $x$ is a nonstandard type-III interior point, 
then it is also a nonstandard type-II interior point. In other words, we have 
$\mbox{int}^{\tiny \mbox{(I)}}(A)\subseteq\mbox{int}(A)$ and 
$\mbox{int}^{\tiny \mbox{(III)}}(A)\subseteq\mbox{int}^{\tiny \mbox{(II)}}(A)$.

\item Suppose that pseudo-metric $d$ satisfies the null equalities and 
$X$ owns the null decomposition. Then 
Proposition~\ref{ch1p273} shows that the concepts of 
nonstandard interior and nonstandard type-I interior are equivalent, and the concepts 
of nonstandard type-II interior and nonstandard type-III interior are 
equivalent. In other words, we have 
$\mbox{int}(A)=\mbox{int}^{\tiny \mbox{(I)}}(A)$ and 
$\mbox{int}^{\tiny \mbox{(II)}}(A)=\mbox{int}^{\tiny \mbox{(III)}}(A)$.
\end{itemize}
}\end{Rem}

\begin{Ex}
{\em
Let ${\cal I}$ be the set of all closed intervals with the null set 
$\Omega =\{[-k,k]:k\geq 0\}$. Given 
$\omega =[-k,k]$ for $k\neq 0$, i.e., $k>0$, we can write $k=k_{1}+k_{2}$ 
with $k_{1},k_{2}>0$. Then, we have 
\[\omega =[-k,k]=[-(k_{1}+k_{2}),k_{1}+k_{2}]=
[-k_{1},k_{1}]\oplus [-k_{2},k_{2}],\]
where $\omega_{1}=[-k_{1},k_{1}],\omega_{2}=[-k_{2},k_{2}]\in\Omega$. 
This says that $\Omega\oplus\Omega =\Omega$ under this space $X$.
Therefore, the assumption in Remark~\ref{ch1r239} (ii) is automatically
satisfied under this space $X$.
}\end{Ex}

\begin{Rem}{\label{ch1r396}}
{\em
Although $B(x;\epsilon )\oplus\Omega\subseteq B(x;\epsilon )$ as shown in 
Proposition~\ref{ch2p91} is satisfied under the assumption of null equalities, 
the set $B(x;\epsilon )\oplus\Omega$ does not necessarily contain the center $x$
unless $x$ has the null decomposition (which will be shown below). 
Therefore, it can happen that there exists an open ball such that 
$B(x;\epsilon )\oplus\Omega$ is contained in $A$ even though the center 
$x$ is not in $A$. In this situation, we do not say that $x$ is a 
nonstandard type-I interior point, since $x$ is not in $A$. 
However, if the center $x$ has the null decomposition, 
then $B(x;\epsilon )\oplus\Omega$ will contain 
the center $x$, since $x=\bar{x}\oplus\omega$ for some $\bar{x}\in X$ 
and $\omega\in\Omega$ satisfying 
\[\epsilon >0=d(x,x)=d(x,\bar{x}\oplus\omega )\geq d(x,\bar{x}),\]
which says that $\bar{x}\in B(x;\epsilon )$, i.e., $x=\bar{x}\oplus\omega
\in B(x;\epsilon )\oplus\Omega$.
}\end{Rem}

Based on Remark~\ref{ch1r396}, we can define the concepts of pseudo-interior point.

\begin{Def}
{\em 
Let $A$ be a subset of a pseudo-metric space $(X,d)$ on a nonstandard vector space $X$. 
\begin{itemize}
\item A point $x_{0}\in X$ is said to be a {\em nonstandard pseudo-interior point} 
of $A$ if there exists $\epsilon >0$ such that $B(x_{0};\epsilon )\subseteq A$. 
The collection of all interior points of $A$ is 
called the {\em nonstandard pseudo-interior} of $A$ and is denoted by 
$\mbox{pint}(A)$.

\item A point $x_{0}\in X$ is said to be a 
{\em nonstandard type-I pseudo-interior point} of $A$ if there exists 
$\epsilon >0$ such that 
$B(x_{0};\epsilon )\oplus\Omega\subseteq A$. The collection of all 
nonstandard type-I pseudo-interior points of $A$ is called the 
{\em nonstandard type-I pseudo-interior} of $A$ and is 
denoted by $\mbox{pint}^{\tiny \mbox{(I)}}(A)$.

\item A point $x_{0}\in X$ is said to be a 
{\em nonstandard type-II pseudo-interior 
point} of $A$ if there exists $\epsilon >0$ such that 
$B(x_{0};\epsilon )\subseteq A\oplus\Omega$. The collection of all 
nonstandard type-II pseudo-interior points of $A$ is called the 
{\em nonstandard type-II pseudo-interior} of $A$ and is 
denoted by $\mbox{pint}^{\tiny \mbox{(II)}}(A)$.

\item A point $x_{0}\in X$ is said to be a 
{\em nonstandard type-III pseudo-interior 
point} of $A$ if there exists $\epsilon >0$ such that 
$B(x_{0};\epsilon )\oplus\Omega\subseteq A\oplus\Omega$. The collection of all 
nonstandard type-III pseudo-interior points of $A$ is called the 
{\em nonstandard type-III pseudo-interior} of $A$ and is 
denoted by $\mbox{pint}^{\tiny \mbox{(III)}}(A)$.
\end{itemize}
}\end{Def}

\begin{Rem}{\label{ch2r501}}
{\em
Let $(X,d)$ be a pseudo-metric space on a nonstandard vector space $X$. 
We have the following observations.
\begin{itemize}
\item Each type of nonstandard pesudo-interior points of $A$ 
does not necessarily belong to $A$. However, each type of nonstandard 
interior points of $A$ is the corresponding nonstandard 
pesudo-interior points of $A$. In other words, we have 
$\mbox{int}^{\tiny \mbox{(I)}}(A)\subseteq\mbox{pint}^{\tiny \mbox{(I)}}(A)$,
$\mbox{int}^{\tiny \mbox{(II)}}(A)\subseteq\mbox{pint}^{\tiny \mbox{(II)}}(A)$ and 
$\mbox{int}^{\tiny \mbox{(III)}}(A)\subseteq\mbox{pint}^{\tiny \mbox{(III)}}(A)$.

\item If the pseudo-metric $d$ satisfies the null inequalities and 
$X$ owns the null decomposition, then,
from Proposition~\ref{ch1p273}, we have $\mbox{pint}^{\tiny \mbox{(I)}}(A)
\subseteq\mbox{int}(A)\subseteq A$, since $x_{0}\in B(x_{0};\epsilon )
\subseteq B(x_{0};\epsilon )\oplus\Omega\subseteq A$.

\item If $A\oplus\Omega\subseteq A$, then we have 
$\mbox{pint}^{\tiny \mbox{(II)}}(A)\subseteq\mbox{int}(A)\subseteq A$.

\item If $A\oplus\Omega\subseteq A$, 
the pseudo-metric $d$ satisfies the null inequalities and 
$X$ owns the null decomposition, then,
from Proposition~\ref{ch1p273}, 
we have $\mbox{pint}^{\tiny \mbox{(III)}}(A)\subseteq\mbox{int}(A)\subseteq A$.
\end{itemize}
}\end{Rem}

\begin{Def}
{\em
Let $A$ be a subset of a pseudo-metric space $(X,d)$ on a nonstandard vector space $X$. 
The set $A$ is said to be {\em nonstandardly open} if $A=\mbox{int}(A)$. 
The set $A$ is said to be {\em nonstandardly type-I-open} if 
$A=\mbox{int}^{\tiny \mbox{(I)}}(A)$. The set $A$ is said to be {\em nonstandardly 
type-II-open} if $A=\mbox{int}^{\tiny \mbox{(II)}}(A)$. The set $A$ is said to be 
{\em nonstandardly type-III-open} if $A=\mbox{int}^{\tiny \mbox{(III)}}(A)$.
We can similarly define the {\em nonstandardly pseudo-open}, 
{\em nonstandardly type-I pseudo-open}, {\em nonstandardly type-II pseudo-open} 
and {\em nonstandardly type-III pseudo-open} set.
}\end{Def}

\begin{Rem}
{\em
Let $A$ be a subset of a pseudo-metric space $(X,d)$ on a nonstandard vector space $X$. 
We have the following observations.
\begin{itemize}
\item From Remark~\ref{ch2r501}, we see that if $A$ is nonstandardly type-I-open, 
then we have $A\subseteq\mbox{pint}^{\tiny \mbox{(I)}}(A)$. 
We can have the similar observations for the other types of openness.

\item It is clear that $\mbox{int}^{\tiny \mbox{(I)}}(A)\subseteq A$. 
Let $O$ be any nonstandardly type-I-open subset of $A$. 
Then $O=\mbox{int}^{\tiny \mbox{(I)}}(O)\subseteq\mbox{int}^{\tiny \mbox{(I)}}(A)$. 
This says that $\mbox{int}^{\tiny \mbox{(I)}}(A)$ is the largest nonstandardly
type-I-open set contained in $A$. Similarly, $\mbox{int}(A)$ is the largest 
open set contained in $A$, $\mbox{int}^{\tiny \mbox{(II)}}(A)$ is the largest nonstandardly
type-II-open set contained in $A$ and $\mbox{int}^{\tiny \mbox{(III)}}(A)$ is the largest 
nonstandardly type-III-open set contained in $A$.
\end{itemize}
}\end{Rem}

For convenience, we adopt $\emptyset\oplus\Omega =\emptyset$.

\begin{Rem}{\label{ch1r430}}
{\em
Let $A$ be a subset of a pseudo-metric space $(X,d)$ on a nonstandard vector space $X$. 
We consider the extreme cases of the empty set $\emptyset$ and whole set $X$.
\begin{enumerate}
\item [(a)] Since the empty set $\emptyset$ contains no elements, it means that 
$\emptyset$ is nonstandardly open and pseudo-open 
(we can regard the empty set as an open ball). 
Obviously, $X$ is also nonstandardly open and pseudo-open, since $x\in B\subseteq X$ 
for any open ball $B$,
i.e., $X\subseteq\mbox{int}(X)$ and $X\subseteq\mbox{pint}(X)$.

\item [(b)] Since $\emptyset\oplus\Omega =\emptyset\subseteq\emptyset$, 
the emptyset $\emptyset$ is nonstandardly type-I-open and type-I-pseudo-open. 
Obviously, $X$ is also nonstandardly type-I-open and type-I-pseudo-open, 
since $x\in B\oplus\Omega\subseteq X$ for any open ball $B$,
i.e., $X\subseteq\mbox{int}^{\tiny \mbox{(I)}}(X)$ and 
$X\subseteq\mbox{pint}^{\tiny \mbox{(I)}}(X)$. 

\item [(c)] Since $\emptyset\subseteq\emptyset =\Omega\oplus\emptyset$, it 
means that $\emptyset$ is nonstandardly type-II-open and type-II-pseudo-open. 
However, $X$ is not nonstandardly type-II-open or type-II-pseudo-open
in general unless $X\subseteq X\oplus\Omega$. If we assume that $X$ owns the 
null decomposition, then $X$ will become a nonstandardly type-II-open and 
and type-II-pseudo-open set.
Indeed, for any $x\in X$ and any open ball $B$, we have 
$x\in B\subseteq X\subseteq X\oplus\Omega$ by the null decomposition,
i.e., $X\subseteq\mbox{int}^{\tiny \mbox{(II)}}(X)$ and $X\subseteq\mbox{pint}^{\tiny \mbox{(II)}}(X)$.

\item [(d)] Since $\emptyset\oplus\Omega\subseteq\Omega\oplus\emptyset$, 
it means that $\emptyset$ is nonstandardly type-III-open and type-III-pseudo-open. 
Now, for any $x\in X$ and any open ball $B$, we have $x\in B\subseteq X$, 
which syas that $B\oplus\Omega\subseteq X\oplus\Omega$, i.e., 
$X\subseteq\mbox{int}^{\tiny \mbox{(III)}}(X)$ and $X\subseteq\mbox{pint}^{\tiny \mbox{(III)}}(X)$.
This shows that $X$ is nonstandardly type-III-open and type-III-pseudo-open.
\end{enumerate}
}\end{Rem}

\begin{Pro}{\label{ch2p2}}
Let $A$ be a nonstandardly pseudo-open, or nonstandardly type-I-pseudo-open, 
or nonstandardly type-II-pseudo-open, or nonstandardly 
type-III-pseudo-open subset of a pseudo-metric 
space $(X,d)$ on a nonstandard vector space $X$. 
Then the following statements hold true.
\begin{enumerate}
\item [{\em (i)}] If the pseudo-metric $d$
satisfies the null inequalities, then $a\in A$ implies $a\oplus\omega\in A$ 
for any $\omega\in\Omega$. If we further assume that 
the metric $d$ satisfies the null equalities. 
Then $a\in A$ if and only if $a\oplus\omega\in A$ 
for any $\omega\in\Omega$.

\item [{\em (ii)}] Suppose that the pseudo-metric $d$ satisfies the null equalities.
Then we have $A\oplus\Omega\subseteq A$
and $A\oplus\omega\subseteq A$
for any $\omega\in\Omega$, and 
$a\oplus\omega\in A\oplus\omega$ implies $a\in A$ for any 
$\omega\in\Omega$.

\item [{\em (iii)}]  Suppose that $X$ owns the null decomposition and 
the pseudo-metric $d$ satisfies the null equalities. Then we have $A=A\oplus\Omega$.
If we further assume that $\Omega$ owns the self-decomposition, then 
we also have $A=A\oplus\omega$ for any $\omega\in\Omega$.
\end{enumerate}
\end{Pro}
\begin{Proof}
To prove part (i), we consider the case of nonstandardly type-III-pseudo-open. 
For $a\in A=\mbox{pint}^{\tiny \mbox{(III)}}(A)$, 
by definition, there exists $\epsilon >0$ such that $B(a;\epsilon )\oplus
\Omega\subseteq A\oplus\Omega$. From Proposition~\ref{ch2p1}, we also have 
$B(a\oplus\omega ;\epsilon )\oplus\Omega\subseteq A\oplus\Omega$, which says 
that $a\oplus\omega\in\mbox{pint}^{\tiny \mbox{(III)}}(A)=A$. 
Suppose that the metric $d$ satisfies the null equalities. If 
$a\oplus\omega\in A=\mbox{pint}^{\tiny \mbox{(III)}}(A)$, there exists $\epsilon >0$ such that 
$B(a\oplus\omega ;\epsilon )\oplus\Omega\subseteq A\oplus\Omega$. 
By Proposition~\ref{ch2p1}, we also see that $a\in\mbox{pint}^{\tiny \mbox{(III)}}(A)=A$.
The similar arguments can also apply to the other cases of openness.

To prove part (ii), we also consider the case of nonstandardly type-III-pseudo-open.
If $a\in A\oplus\omega$, then $a=\hat{a}\oplus\omega$ for some 
$\hat{a}\in A=\mbox{pint}^{\tiny \mbox{(III)}}(A)$. Therefore, there exists $\epsilon >0$ such 
that $B(\hat{a};\epsilon )\oplus\Omega\subseteq A\oplus\Omega$. Since 
$B(a;\epsilon )=B(a\oplus\omega ;\epsilon )=B(\hat{a};\epsilon )$ by 
Proposition~\ref{ch2p1}, we see that $B(a;\epsilon )\oplus\Omega\subseteq 
A\oplus\Omega$, i.e., $a\in\mbox{pint}^{\tiny \mbox{(III)}}(A)=A$. Now, for $a\in A\oplus\Omega$,
we see that $a\in A\oplus\omega$ for some $\omega\in A$, which implies
$a\in A$. Therefore, we obtain $A\oplus\Omega\subseteq A$.
Now, for $a\oplus\omega\in A\oplus\omega\subseteq A\oplus\Omega\subseteq A$, 
we also have $a\oplus\omega\in A$. From (i), we obtain $a\in A$.

To prove part (iii), for $a\in A$, since $X$ owns the 
null decomposition, we have $a=\hat{a}\oplus\hat{\omega}_{0}$ for some 
$\hat{a}\in X$ and $\hat{\omega}_{0}\in\Omega$.
Since $a\in A=\mbox{pint}^{\tiny \mbox{(III)}}(A)$, there exists $\epsilon >0$ 
such that $B(\hat{a}\oplus\hat{\omega}_{0};\epsilon )\oplus\Omega
\subseteq A\oplus\Omega$. 
From Proposition~\ref{ch2p1}, we also have $B(\hat{a};\epsilon )\oplus\Omega
=B(\hat{a}\oplus\hat{\omega}_{0};\epsilon )\oplus\Omega\subseteq A\oplus\Omega$, i.e., 
$\hat{a}\in\mbox{pint}^{\tiny \mbox{(III)}}(A)=A$. This shows that $a=\hat{a}\oplus\hat{\omega}_{0}
\in A\oplus\Omega$. Therefore, from (ii), we conclude that $A=A\oplus\Omega$.
We furthre assume that $\Omega$ owns the self-decomposition. Then, given any 
$\omega\in\Omega$, we have
$\hat{\omega}_{0}=\omega\oplus\hat{\omega}_{1}$ for some $\hat{\omega}_{1}
\in\Omega$. From Proposition~\ref{ch2p1} again, we have 
$B(\hat{a}\oplus\hat{\omega}_{1};\epsilon )\oplus\Omega=
B(\hat{a};\epsilon )\oplus\Omega\subseteq A\oplus\Omega$, i.e., 
$\hat{a}\oplus\hat{\omega}_{1}\in\mbox{pint}^{\tiny \mbox{(III)}}(A)=A$. Therefore, we 
obtain $a=\hat{a}\oplus\hat{\omega}_{0}=
(\hat{a}\oplus\hat{\omega}_{1})\oplus\omega\in A\oplus\omega$. From (ii),
we conclude that $A=A\oplus\omega$ for any $\omega$. This completes the proof.
\end{Proof}

\begin{Pro}{\label{ch1r405}}
Let $(X,d)$ be a pseudo-metric space on a nonstandard vector space $X$. 
Then the following statements hold true.
\begin{enumerate}
\item [{\em (i)}] We have $\mbox{\em int}^{\tiny \mbox{\em (I)}}(A)\oplus\Omega\subseteq A$. Moreover, 
if $A$ is a nonstandardly type-I open subset of $X$, then $A\oplus\Omega
\subseteq A$.

\item [{\em (ii)}] We have $\mbox{\em int}^{\tiny \mbox{\em (II)}}(A)\subseteq A\oplus\Omega$. Moreover, 
if $A$ is a nonstandardly type-II-open subset of $X$, 
then $A\subseteq A\oplus\Omega$.

\item [{\em (iii)}] Suppose that the pseudo-metric $d$ satisfies the null equalities.
We have $(\mbox{\em int}^{\tiny \mbox{\em (II)}}(A))^{c}\oplus\Omega\subseteq (\mbox{\em int}^{\tiny \mbox{\em (II)}}(A))^{c}$.
Moreover, if $A$ is nonstandardly type-II-open subset of $X$, 
then $A^{c}\oplus\Omega\subseteq A^{c}$.

\item [{\em (iv)}] Suppose that the pseudo-metric $d$ satisfies the null equalities.
We have $\mbox{\em int}(A)\oplus\Omega\subseteq A$ .
Moreover, if $A$ is a nonstandardly open subset of $X$, 
then $A\oplus\Omega\subseteq A$.

\item [{\em (v)}] Suppose that the pseudo-metric $d$ satisfies the null inequalities
and $X$ owns the null decomposition. 
We have $\mbox{\em int}^{\tiny \mbox{\em (III)}}(A)\subseteq A\oplus\Omega$. Moreover, 
if $A$ is a nonstandardly type-III-open subset of $X$, then 
$A\subseteq A\oplus\Omega$.

\item [{\em (vi)}] Suppose that the pseudo-metric $d$ satisfies the null equalities.
We have $(\mbox{\em int}^{\tiny \mbox{\em (III)}}(A))^{c}\oplus\Omega\subseteq (\mbox{\em int}^{\tiny \mbox{\em (III)}}(A))^{c}$.
Moreover, if $A$ is nonstandardly type-III-open subset of $X$, 
then $A^{c}\oplus\Omega\subseteq A^{c}$.
\end{enumerate}
\end{Pro}
\begin{Proof}
To prove part (i), it suffices to prove the case of null set $\Omega$.
For any $x\in\mbox{int}^{\tiny \mbox{(I)}}(A)$, there exists an open ball $B(x;\epsilon )$
such that $B(x;\epsilon )\oplus\Omega\subseteq A$. Since $x\in B(x;\epsilon )$,
we have $x\oplus\Omega\subseteq B(x;\epsilon )\oplus\Omega\subseteq A$.
This shows that $\mbox{int}^{\tiny \mbox{(I)}}(A)\oplus\Omega\subseteq A$. 

To prove part (ii), for any $x\in\mbox{int}^{\tiny \mbox{(II)}}(A)$, 
there exists an open ball 
$B(x;\epsilon )$ such that $B(x;\epsilon )\subseteq A\oplus\Omega$. 
Then, we have $x\in A\oplus\Omega$, since $x\in B(x;\epsilon )$.
This shows that $\mbox{int}^{\tiny \mbox{(II)}}(A)\subseteq A\oplus\Omega$.

To prove part (iii), for any $x\in (\mbox{int}^{\tiny \mbox{(II)}}(A))^{c}\oplus\Omega$, we have 
$x=\hat{x}\oplus\hat{\omega}$ for some $\hat{x}\in (\mbox{int}^{\tiny \mbox{(II)}}(A))^{c}$
and $\hat{\omega}\in\Omega$. By definition, we see that $B(\hat{x};\epsilon )
\not\subseteq A\oplus\Omega$ for every $\epsilon >0$. Since $B(x;\epsilon )=
B(\hat{x}\oplus\hat{\omega};\epsilon )=B(\hat{x};\epsilon )$ by 
Proposition~\ref{ch2p1}, we also have $B(x;\epsilon )
\not\subseteq A\oplus\Omega$ for every $\epsilon >0$. This says that $x$ is 
not a nonstandardly type-II interior point of $A$, i.e., $x\not\in 
\mbox{int}^{\tiny \mbox{(II)}}(A)=A$, which shows that $x\in A^{c}$.

To prove part (iv), for any $x\in\mbox{int}(A)$, there exists an open ball 
$B(x;\epsilon )$ contained in $A$. From Proposition~\ref{ch2p91}, we have 
$x\oplus\Omega\subseteq B(x;\epsilon )\oplus\Omega\subseteq B(x;\epsilon )
\subseteq A$. This shows that $\mbox{int}(A)\oplus\Omega\subseteq A$. 

To prove part (v), for any $x\in\mbox{int}^{\tiny \mbox{(III)}}(A)$,
 there exists an open ball 
$B(x;\epsilon )$ such that $B(x;\epsilon )\oplus\Omega\subseteq A\oplus\Omega$. 
From Proposition~\ref{ch1p273}, we have $x\in B(x;\epsilon )\subseteq 
B(x;\epsilon )\oplus\Omega\subseteq A\oplus\Omega$.
This shows that $\mbox{int}^{\tiny \mbox{(III)}}(A)\subseteq A\oplus\Omega$.

To prove part (vi), for any $x\in (\mbox{int}^{\tiny \mbox{(III)}}(A))^{c}\oplus\Omega$, we have 
$x=\hat{x}\oplus\hat{\omega}$ for some $\hat{x}\in (\mbox{int}^{\tiny \mbox{(III)}}(A))^{c}$
and $\hat{\omega}\in\Omega$. By the arguments of (iii), we see that $B(x;\epsilon )
\oplus\Omega =
B(\hat{x};\epsilon )\oplus\Omega\not\subseteq A\oplus\Omega$ for every $\epsilon >0$. 
This says that $x\not\in\mbox{int}^{\tiny \mbox{(III)}}(A)=A$, which shows that $x\in A^{c}$.
We complete the proof.
\end{Proof}

\begin{Pro}{\label{ch1p409}}
Let $(X,d)$ be a pseudo-metric space on a nonstandard vector space $X$.
Then, the following statements hold true.
\begin{enumerate}
\item [{\em (i)}] Suppose that the pseudo-metric $d$ satisfies the null inequalities. 
If $A$ is nonstandardly open, then $A$ is nonstandardly type-III-open.

\item [{\em (ii)}] Suppose that the metric $d$ satisfies the null equalities.
If $A$ is nonstandardly open, then $A$ is nonstandardly type-I-open; 
if $A$ is nonstandardly type-II-open, then $A$ is also 
nonstandardly type-III-open.

\item [{\em (iii)}] Suppose that the pseudo-metric $d$ satisfies the null inequalities 
and $\Omega\oplus\Omega =\Omega$. If $A$ is 
nonstandardly type-I-open, then $A$ is nonstandardly type-III-open.

\item [{\em (iv)}] Suppose that the pseudo-metric $d$ satisfies the null equalities.
If $A$ is simultaneously nonstandardly type-I and type-II-open, 
then $A$ is simultaneously nonstandardly open and nonstandardly type-III-open.

\item [{\em (v)}] Suppose that $X$ own the null decomposition and 
the pseudo-metric $d$ satisfies the null equalities.
If $A$ is simultaneously nonstandardly open and nonstandardly type-III-open, 
then $A$ is also nonstandardly type-I-open.

\item [{\em (vi)}] Suppose that $\Omega\oplus\Omega =\Omega$ and the pseudo-metric 
$d$ satisfies null equalities. If $A$ is nonstandardly open, then 
$A$ is nonstandardly type-I-open if and only if $A$ is nonstandardly 
type-III-open.

\item [{\em (vii)}] Suppose that the space $X$ owns the null decomposition
and the pseudo-metric $d$ satisfies the null equalities. 
Then $A$ is nonstandardly open if and only if $A$ is 
nonstandard type-I-open, and $A$ is nonstandard type-II-open if and only if
$A$ is nonstandard type-III-open.

\item [{\em (viii)}] Suppose that the space $X$ owns the null decomposition, 
the pseudo-metric $d$ satisfies the null equalities and $\Omega\oplus\Omega =\Omega$. 
If $A$ is nonstandardly open or nonstandardly type-I open, 
then $A$ is simultaneously open, nonstandardly type-I-open, type-II-open and 
type-III-open.
\end{enumerate}
\end{Pro}
\begin{Proof}
To prove part (i), from Remark~\ref{ch1r239} (i), we see that 
$\mbox{int}(A)\subseteq\mbox{int}^{\tiny \mbox{(III)}}(A)$.
If $A$ is nonstandardly open, then 
$A=\mbox{int}(A)\subseteq\mbox{int}^{\tiny \mbox{(III)}}(A)$, i.e., 
$A=\mbox{int}^{\tiny \mbox{(III)}}(A)$. This shows that $A$ is 
nonstandardly type-III open.
Part (ii) follows from Remark~\ref{ch1r239} (iii) immediately.
Part (iii) follows from Remark~\ref{ch1r239} (ii) immediately.
To prove part (iv), if $A$ is simultaneously nonstandardly type-I-open and type-II-open, 
then, by Proposition~\ref{ch1r405} (i) and (ii), we have $B(x_{0};\epsilon )
\subseteq A\oplus\Omega =A$ and $B(x_{0};\epsilon )\oplus\Omega\subseteq A=
A\oplus\Omega$. This proves the result.
To prove part (v), if $A$ is nonstandardly open and nonstandardly type-III-open, 
then $B(x_{0};\epsilon )\oplus\Omega\subseteq A\oplus\Omega\subseteq A$ by 
Proposition~\ref{ch1r405} (iii) and (iv).
Part (vi) follows from (iii) and (v) immediately. 
Part (vii) follows from  Remark~\ref{ch1r239} (v) immediately. 
Part (viii) follows from (vi) and (vii) immediately.
This completes the proof.
\end{Proof}

\begin{Pro}{\label{ch1p410}}
Let $X$ be a pseudo-metric space on a nonstandard vector space $X$ such that the 
pseudo-metric $d$ satisfies the null inequalities. 
Suppose that the space $X$ owns the null decomposition, 
the pseudo-metric $d$ satisfies the null equalities and $\Omega\oplus\Omega =\Omega$. 
If $A$ is nonstandardly open, then $A=A\oplus\Omega$.
\end{Pro}
\begin{Proof}
The result follows from Propositions~\ref{ch1p409} (viii) and \ref{ch1r405} (iii)
immediately.
\end{Proof}

\begin{Pro}{\label{ch1p255}}
Let $(X,d)$ be a pseudo-metric space on a nonstandard vector space $X$ such that the 
pseudo-metric $d$ satisfies the null inequalities, 
and let $A$ be a subset of $X$. Then we have 
$\mbox{\em int}(\mbox{\em int}(A))=\mbox{\em int}(A)$ and 
$\mbox{\em int}^{\tiny \mbox{\em (I)}}(\mbox{\em int}^{\tiny \mbox{\em (I)}}(A))=\mbox{\em int}^{\tiny \mbox{\em (I)}}(A)$.
In other words, $\mbox{\em int}(A)$ is nonstandardly open and $\mbox{\em int}^{\tiny \mbox{\em (I)}}(A)$
is nonstandardly type-I open. If we further assume that the pseudo-metric $d$ 
satisfies the null equalities, then we have
$\mbox{\em int}^{\tiny \mbox{\em (II)}}(\mbox{\em int}^{\tiny \mbox{\em (II)}}(A))=\mbox{\em int}^{\tiny \mbox{\em (II)}}(A)$
and $\mbox{\em int}^{\tiny \mbox{\em (III)}}(\mbox{\em int}^{\tiny \mbox{\em (III)}}(A))=\mbox{\em int}^{\tiny \mbox{\em (III)}}(A)$.
In other words, $\mbox{\em int}^{\tiny \mbox{\em (II)}}(A)$ is nonstandardly type-II open and 
$\mbox{\em int}^{\tiny \mbox{\em (III)}}(A)$ is nonstandardly type-III open.
\end{Pro}
\begin{Proof}
We consider the case of nonstandardly type-I-open set.
It will be enough to show the inclusion $\mbox{int}^{\tiny \mbox{(I)}}(A)\subseteq
\mbox{int}^{\tiny \mbox{(I)}}(\mbox{int}^{\tiny \mbox{(I)}}(A))$. Suppose that $x\in\mbox{int}^{\tiny \mbox{(I)}}(A)$. Then 
there exist $\epsilon >0$ such that $B(x;\epsilon )\oplus\Omega\subseteq A$. 
We want to claim that each element of $B(x;\epsilon )\oplus\Omega$ is a 
nonstandard type-I interior point of $A$. We take the element 
$\hat{y}=y\oplus\omega\in B(x;\epsilon )\oplus\Omega$, where $y\in 
B(x;\epsilon )$ and $\omega\in\Omega$, i.e., $d(x,y)<\epsilon$. 
Let $\hat{\epsilon}=d(x,y)$. Then we consider the open ball
$B(\hat{y},\epsilon -\hat{\epsilon})$ centered at $\hat{y}$ with radius
$\epsilon -\hat{\epsilon}$. For $z\in B(\hat{y},\epsilon -\hat{\epsilon})$,
from condition (iv) in Definition~\ref{ch1d7}, we have 
\begin{equation}{\label{ch1eq255}}
d(y,z)\leq d(y\oplus\omega ,z)=d(\hat{y},z)<\epsilon -\hat{\epsilon}.
\end{equation}
Therefore, we obtain
\begin{equation}{\label{ch1eq256}}
d(x,z)\leq d(x,y)+d(y,z)=\hat{\epsilon}+d(y,z)<\hat{\epsilon}+
\epsilon -\hat{\epsilon}=\epsilon .
\end{equation}
This shows that $z\in B(x,\epsilon )$, i.e., $B(\hat{y},\epsilon 
-\hat{\epsilon})\subseteq B(x,\epsilon )$. Then we have 
\[B(\hat{y},\epsilon -\hat{\epsilon})\oplus\Omega
\subseteq B(x;\epsilon )\oplus\Omega\subseteq A,\]
which says that $\hat{y}\in\mbox{int}^{\tiny \mbox{(I)}}(A)$. This shows that 
$B(x;\epsilon )\oplus\Omega\subseteq\mbox{int}^{\tiny \mbox{(I)}}(A)$, i.e., 
$x\in\mbox{int}^{\tiny \mbox{(I)}}(\mbox{int}^{\tiny \mbox{(I)}}(A))$. Without considering the null 
set $\Omega$, we can also use the above similar arguments to show 
$\mbox{int}(A)\subseteq \mbox{int}(\mbox{int}(A))$. 

For the case of nonstandardly type-II-open set, if $x\in\mbox{int}^{\tiny \mbox{(II)}}(A)$, 
there exists $\epsilon >0$ such that 
$B(x;\epsilon )\subseteq A\oplus\Omega$. Let $\hat{y}\in B(x;\epsilon )$ and
$\hat{\epsilon}=d(x,\hat{y})$. Therefore, we have $\hat{y}=a\oplus\omega$ 
for some $a\in A$ and $\omega\in\Omega$. We consider the open ball
$B(\hat{y},\epsilon -\hat{\epsilon})$. 
For $z\in B(\hat{y},\epsilon -\hat{\epsilon})$, we have 
\[d(x,z)\leq d(x,\hat{y})+d(\hat{y},z)=\hat{\epsilon}+d(\hat{y},z)
<\hat{\epsilon}+\epsilon -\hat{\epsilon}=\epsilon .\]
This shows that $z\in B(x,\epsilon )$, i.e., $B(\hat{y},\epsilon 
-\hat{\epsilon})\subseteq B(x,\epsilon )\subseteq A\oplus\Omega$.
Now we want to claim that $B(a;\epsilon -\hat{\epsilon})
\subseteq B(\hat{y};\epsilon -\hat{\epsilon})$.
Suppose that $\hat{x}\in B(a;\epsilon -\hat{\epsilon})$. Since
the metric $d$ satisfies thenull equalities, we have 
\[d(\hat{x},\hat{y})=d(\hat{x},a\oplus\omega )=d(\hat{x},a)<\epsilon -
\hat{\epsilon},\]
which says that $\hat{x}\in B(\hat{y};\epsilon -\hat{\epsilon})$. Therefore, 
we obtain $B(a;\epsilon -\hat{\epsilon})\subseteq 
B(\hat{y};\epsilon -\hat{\epsilon})\subseteq A\oplus\Omega$. This shows
that $a\in\mbox{int}^{\tiny \mbox{(II)}}(A)$, i.e., 
$\hat{y}=a\oplus\omega\in\mbox{int}^{\tiny \mbox{(II)}}(A)\oplus\Omega$. Therefore, we obtain
$B(x;\epsilon )\subseteq\mbox{int}^{\tiny \mbox{(II)}}(A)\oplus\Omega$, which implies that 
$x\in\mbox{int}^{\tiny \mbox{(II)}}(\mbox{int}^{\tiny \mbox{(II)}}(A))$. Therefore, we obtain 
the inclusion $\mbox{int}^{\tiny \mbox{(II)}}(A)\subseteq
\mbox{int}^{\tiny \mbox{(I)}}(\mbox{int}^{\tiny \mbox{(II)}}(A))$.

Finally, for the case of nonstandardly type-III-open set, 
if $x\in\mbox{int}^{\tiny \mbox{(III)}}(A)$, there exists $\epsilon >0$ 
such that $B(x;\epsilon )\oplus\Omega\subseteq A\oplus\Omega$. 
For $\hat{y}\in B(x;\epsilon )$ and $\omega\in\Omega$, 
we have $\hat{y}\oplus\omega=a\oplus\hat{\omega}$ for some 
$a\in A$ and $\hat{\omega}\in\Omega$. Let $\hat{\epsilon}=d(x,\hat{y})$. 
From the above arguments, we have $B(\hat{y},\epsilon -\hat{\epsilon})
\subseteq B(x;\epsilon )$, i.e., 
\[B(\hat{y},\epsilon -\hat{\epsilon})\oplus\Omega
\subseteq B(x;\epsilon )\oplus\Omega\subseteq A\oplus\Omega .\]
We want to claim that 
$B(a;\epsilon -\hat{\epsilon})\subseteq B(\hat{y};\epsilon -\hat{\epsilon})$.
Suppose that $\hat{x}\in B(a;\epsilon -\hat{\epsilon})$. Since
the metric $d$ satisfies the null equalities, we have 
\[d(\hat{x},\hat{y})=d(\hat{x},\hat{y}\oplus\omega )=
d(\hat{x},a\oplus\hat{\omega})=d(\hat{x},a)<\epsilon -\hat{\epsilon},\]
which says that $\hat{x}\in B(\hat{y};\epsilon -\hat{\epsilon})$.
Therefore, we obtain 
\[B(a;\epsilon -\hat{\epsilon})\oplus\Omega\subseteq 
B(\hat{y};\epsilon -\hat{\epsilon})\oplus\Omega\subseteq A\oplus\Omega .\]
This shows that $a\in\mbox{int}^{\tiny \mbox{(III)}}(A)$, i.e., we have 
$\hat{y}\oplus\omega =a\oplus\hat{\omega}\in\mbox{int}^{\tiny \mbox{(III)}}(A)\oplus\Omega$. 
Therefore, we obtain
$B(x;\epsilon )\oplus\Omega\subseteq\mbox{int}^{\tiny \mbox{(III)}}(A)\oplus\Omega$, which 
implies that $x\in\mbox{int}^{\tiny \mbox{(III)}}(\mbox{int}^{\tiny \mbox{(III)}}(A))$. Therefore, we obtain 
the inclusion $\mbox{int}^{\tiny \mbox{(III)}}(A)\subseteq
\mbox{int}^{\tiny \mbox{(I)}}(\mbox{int}^{\tiny \mbox{(III)}}(A))$. This completes the proof.
\end{Proof}

\begin{Pro}{\label{ch1p267}}
Let $(X,d)$ be a pseudo-metric space on a nonstandard vector space $X$ such that the 
pseudo-metric $d$ satisfies the null inequalities, 
and let $B(x_{0};\epsilon )$ be any open ball centered at $x_{0}$ with radius 
$\epsilon$. Then the following statements hold true.
\begin{enumerate}
\item [{\em (i)}] The open ball $B(x_{0};\epsilon )$ is also an nonstandardly open and 
nonstandardly type-III-open subset of $X$. The result also holds true if 
$(X,d)$ is taken as a pseudo-metric space.

\item [{\em (ii)}] Suppose that the metric $d$ satisfies the null equalities. Then 
the open ball $B(x_{0};\epsilon )$ is also a
nonstandardly type-I-open subset of $X$. 

\item [{\em (iii)}] Suppose that $X$ owns the null decomposition. 
Then the open ball $B(x_{0};\epsilon )$ is 
also a nonstandardly type-II-open subset of $X$.
\end{enumerate}
\end{Pro}
\begin{Proof}
To prove part (i), for any $x\in B(x_{0};\epsilon )$, we have $d(x,x_{0})<\epsilon$.
Let $\hat{\epsilon}=d(x,x_{0})$. We consider the open ball
$B(x;\epsilon -\hat{\epsilon})$ centered at $x$ with radius $\epsilon -
\hat{\epsilon}$. Then, for any $\hat{x}\in B(x;\epsilon -\hat{\epsilon})$,
we have $d(\hat{x},x)<\epsilon -\hat{\epsilon}$. Then we have
\[d(\hat{x},x_{0})\leq d(\hat{x},x)+d(x,x_{0})
=\hat{\epsilon}+d(\hat{x},x)<\hat{\epsilon}+\epsilon -\hat{\epsilon}
=\epsilon ,\]
which means that $\hat{x}\in B(x_{0};\epsilon )$, i.e., 
\begin{equation}{\label{ch1eq412}}
B(x;\epsilon -\hat{\epsilon})\subseteq B(x_{0};\epsilon ).
\end{equation}
This shows that the open ball $B(x_{0};\epsilon )$ is nonstandardly open. Moreover, we also have
$B(x;\epsilon -\hat{\epsilon})\oplus\Omega\subseteq B(x_{0};\epsilon )\oplus\Omega$
from (\ref{ch1eq412}). This says that $B(x_{0};\epsilon )$ is nonstandardly 
type-III-open.

To prove part (ii), suppose that the metric $d$ satisfies the null equalities. 
Then, from Proposition~\ref{ch2p91} and (\ref{ch1eq412}), we have 
$B(x;\epsilon -\hat{\epsilon})\oplus\Omega\subseteq B(x_{0};\epsilon )
\oplus\Omega\subseteq B(x_{0};\epsilon )$. This says that 
$B(x_{0};\epsilon )$ is nonstandardly type-I-open. 

To prove part (iii), suppose that $X$ owns the null decomposition.
From (\ref{ch1eq412}) and Proposition~\ref{ch1p273}, we have 
$B(x;\epsilon -\hat{\epsilon})\subseteq B(x_{0};\epsilon )\subseteq 
B(x_{0};\epsilon )\oplus\Omega$. This says that the open ball 
$B(x_{0};\epsilon )$ is also nonstandardly type-II-open. We complete the proof.
\end{Proof}

\section{Nonstandardly Closed Sets}

In the (conventional) metric space $(X,d)$, let $A$ be a subset of $X$.
A point $x_{0}\in X$ is said to be a closure point of $A$ (or limit point of $A$) 
if every open ball $B(x_{0};\epsilon )$ centered at $x_{0}$ contains points of $A$. 
Equivalently, since $x_{0}\in B(x_{0};\epsilon )$, the point $x_{0}\in X$ is a 
closure point of $A$ if $x_{0}\in A$ or, for $x_{0}\in A^{c}$, every open 
ball $B(x_{0};\epsilon )$ centered at $x_{0}$ contains points of $A$. 
For the pseudo-metric space $(X,d)$ on a nonstandard vector space $X$,
Now different types of closure points are defined below.

\begin{Def}{\label{ch2d507}}
{\em
Let $(X,d)$ be a pseudo-metric space on a nonstandard vector space $X$, 
and let $A$ be a subset of $X$. 
\begin{itemize}
\item A point $x_{0}\in X$ is said to be a {\em nonstandard closure point} of $A$ 
(or {\em nonstandard limit point} of $A$) if every open 
ball $B(x_{0};\epsilon )$ centered at $x_{0}$ contains points of $A$. 
The collection of all closure points of $A$ is 
called the {\em nonstandard closure} of $A$ and is denoted by $\mbox{cl}(A)$

\item A point $x_{0}\in X$ is said to be a {\em nonstandard type-I closure 
point} of $A$ (or {\em nonstandard type-I limit point} of $A$) if $x_{0}\in A$, or if
$x_{0}\in A^{c}$ and, for every open ball $B(x_{0};\epsilon )$ centered at $x_{0}$, 
the set $B(x_{0};\epsilon )\oplus\Omega$ contains points of $A$. The collection 
of all nonstandard type-I closure points of $A$ is called the 
{\em nonstandard type-I closure} of $A$ and is denoted by $\mbox{cl}^{\tiny \mbox{(I)}}(A)$.

\item A point $x_{0}\in X$ is said to be a {\em nonstandard type-II closure 
point} of $A$ (or {\em nonstandard type-II limit point} of $A$) if $x_{0}\in A$, or if 
$x_{0}\in A^{c}$ and every open ball $B(x_{0};\epsilon )$ centered at $x_{0}$
contains points of $A\oplus\Omega$. The collection 
of all nonstandard type-II closure points of $A$ is called the 
{\em nonstandard type-II closure} of $A$ and is denoted by $\mbox{cl}^{\tiny \mbox{(II)}}(A)$.

\item A point $x_{0}\in X$ is said to be a {\em nonstandard type-III closure 
point} of $A$ (or {\em nonstandard type-III limit point} of $A$) if $x_{0}\in A$, or 
if $x_{0}\in A^{c}$ and, for every open 
ball $B(x_{0};\epsilon )$ centered at $x_{0}$, the set $B(x_{0};\epsilon )
\oplus\Omega$ contains points of $A\oplus\Omega$. The collection 
of all nonstandard type-III closure points of $A$ is called the 
{\em nonstandard type-III closure} of $A$ and is denoted by $\mbox{cl}^{\tiny \mbox{(III)}}(A)$.
\end{itemize}
}\end{Def}

We also remark that if $X$ happens to be a (conventional) vector space, 
then the four concepts of nonstandard closure point coincide with the 
conventional definition of closure point. The concepts proposed in 
Definition~\ref{ch2d507} are based on the following observations: 
\begin{itemize}
\item For any $x\in A$,
the set $B(x;\epsilon )\oplus\Omega$ does not necessarily contains the 
points of $A$ in the sense of nonstandard type-I closure point, 
unless $B(x;\epsilon )\subseteq B(x;\epsilon )\oplus\Omega$. 
In this case, we have $x\in B(x;\epsilon )$, i.e., $B(x;\epsilon )$
contains points of $A$, which implies $B(x;\epsilon )\oplus\Omega$
contains points of $A$. However, if $x\in A$, 
we still define $x$ as a nonstandard type-I closure point of $A$
as in the conventional case.

\item For any $x\in A$,
the open ball $B(x;\epsilon )$ does not necessarily contains the 
points of $A\oplus\Omega$ in the sense of nonstandard type-II closure point, 
unless $A\subseteq A\oplus\Omega$.
In this case, we have $x\in B(x;\epsilon )$, i.e., $B(x;\epsilon )$
contains points of $A$, which implies $B(x;\epsilon )$
contains points of $A\oplus\Omega$.
However, if $x\in A$, we still define $x$ as a 
nonstandard type-II closure point of $A$
as in the conventional case.

\item For any $x\in A$,
the open ball $B(x;\epsilon )\oplus\Omega$ does not necessarily contains the 
points of $A\oplus\Omega$ in the sense of nonstandard type-III closure point, 
unless $B(x;\epsilon )\subseteq B(x;\epsilon )\oplus\Omega$ and 
$A\subseteq A\oplus\Omega$.
In this case, we have $x\in B(x;\epsilon )$, i.e., $B(x;\epsilon )$
contains points of $A$, which implies $B(x;\epsilon )\oplus\Omega$
contains points of $A\oplus\Omega$.
However, if $x\in A$, we still define $x$ as a 
nonstandard type-III closure point of $A$
as in the conventional case.
\end{itemize}

\begin{Rem}
{\em
Let $(X,d)$ be a pseudo-metric space on a nonstandard vector space $X$ such that the 
pseudo-metric $d$ satisfies the null inequalities. 
We have the following observations
\begin{itemize}
\item We further assume that the pseudo-metric $d$ satisfies the null equalities. 
Proposition~\ref{ch2p91} shows 
the inclusion $B(x_{0};\epsilon )\oplus\Omega\subseteq B(x_{0};\epsilon )$.
This also says that if $x$ is a nonstandard type-I closure point, 
then it is also a closure point, and if $x$ is a nonstandard type-III 
closure point, then it is 
also a nonstandard type-II closure point. 
Now we further assume that $\Omega$ is closed under 
the vector addition, i.e., $\Omega\oplus\Omega\subseteq\Omega$. If $x$ is a
nonstandard type-II closure point, then $B(x;\epsilon )$ contains 
points of $A\oplus\Omega$. Equivalently, 
by adding $\Omega$, we see that $B(x;\epsilon )\oplus\Omega$ contains points of 
$A\oplus\Omega\oplus\Omega\subseteq A\oplus\Omega$. We conclude that 
$B(x;\epsilon )\oplus\Omega$
contains points of $A\oplus\Omega$, i.e.,
$x$ is a nonstandard type-III closure point. This shows that 
the concepts of nonstandard type-II closure point and nonstandard 
type-III closure point are equivalent under the assumption that 
$\Omega$ is closed under the vector addition and the metric 
$d$ satisfies the null equalities.

\item Suppose that $X$ owns the null decomposition.
By Proposition~\ref{ch1p273}, we have the inclusion 
$B(x_{0};\epsilon )\subseteq B(x_{0};\epsilon )\oplus\Omega$.
This also says that if $x$ is a closure point, 
then it is also a nonstandard type-I closure point, and if $x$ is a 
nonstandard type-II closure point, then it is 
also a nonstandard type-III closure point. 

\item Suppose that $X$ owns the null decomposition
and the metric $d$ satisfies the null equalities.
By Proposition~\ref{ch1p273}, the concepts of closure point 
and nonstandard type-I closure point are equivalent, and 
the concepts of nonstandard type-II closure point 
and nonstandard type-III closure point are equivalent.
\end{itemize}
}\end{Rem}

\begin{Rem}{\label{ch1r421}}
{\em
Remark~\ref{ch1r396} says that if the center $x$ has the null decomposition
and the pseudo-metric $d$ satisfies the null equalities, 
then $B(x;\epsilon )\oplus\Omega$
will contain the center $x$. Therefore, if $X$ owns the null decomposition 
and the metric $d$ satisfies the null equalities, 
then we can simply say that a point $x_{0}\in X$ is a nonstandard type-I closure 
point of $A$ if, for every open ball $B(x_{0};\epsilon )$ centered at $x_{0}$, 
the set $B(x_{0};\epsilon )\oplus\Omega$ contains points of $A$.
}\end{Rem}

\begin{Def}
{\em
Let $A$ be a subset of a pseudo-metric space $(X,d)$ on a nonstandard vector space $X$. 
The set $A$ is said to be {\em nonstandardly closed} if and only if $A=\mbox{cl}(A)$.
The set $A$ is said to be {\em nonstandardly type-I-closed} if and only if 
$A=\mbox{cl}^{\tiny \mbox{(I)}}(A)$. The set $A$ is said to be 
{\em nonstandardly type-II-closed} if and only if $A=\mbox{cl}^{\tiny \mbox{(II)}}(A)$.
The set $A$ is said to be {\em nonstandardly type-III-closed} if and only if 
$A=\mbox{cl}^{\tiny \mbox{(III)}}(A)$.
}\end{Def}

\begin{Rem}{\label{ch1r420}}
{\em
We have the following observations.
\begin{enumerate}
\item [(i)] It is clear that $A\subseteq\mbox{cl}(A)$, 
$A\subseteq\mbox{cl}^{\tiny \mbox{(I)}}(A)$, 
$A\subseteq\mbox{cl}^{\tiny \mbox{(II)}}(A)$ and 
$A\subseteq\mbox{cl}^{\tiny \mbox{(III)}}(A)$.

\item [(ii)] It is also clear that the empty set and the whole space 
$X$ are nonstandardly closed, nonstandardly type-I-closed, type-II-closed, and 
type-III-closed. 

\item [(iii)] Let $C$ be any nonstandardly type-I-closed set containing $A$. Then 
$\mbox{cl}^{\tiny \mbox{(I)}}(A)\subseteq\mbox{cl}^{\tiny \mbox{(I)}}(C)=C$. 
This says that $\mbox{cl}^{\tiny \mbox{(I)}}(A)$ 
is the smallest nonstandardly type-I-closed set containing $A$. Similarly, 
$\mbox{cl}(A)$ is the smallest nonstandardly closed set containing $A$, 
$\mbox{cl}^{\tiny \mbox{(II)}}(A)$ 
is the smallest nonstandardly type-II-closed set containing $A$ and 
$\mbox{cl}^{\tiny \mbox{(III)}}(A)$ is the smallest nonstandardly type-III-closed set 
containing $A$.

\item [(iv)] Let $A$ be a subset of $X$ such that, for every $x\in A^{c}$, $x$ is 
not a (resp. nonstandard type-I, type-II and type-III) closure point
of $A$. In other words, by definition,  
if $x\in\mbox{cl}(A)$ (resp. $x\in\mbox{cl}^{\tiny \mbox{(I)}}(A)$, 
$x\in\mbox{cl}^{\tiny \mbox{(II)}}(A)$ and $x\in\mbox{cl}^{\tiny \mbox{(III)}}(A)$) 
then $x\in A$, i.e., $\mbox{cl}(A)\subseteq A$.
This says that $A$ is nonstandardly closed (resp. nonstandardly 
type-I-closed, type-II-closed and type-III-closed) in $X$.
\end{enumerate}
}\end{Rem}

\begin{Rem}
{\em
We have the following observations.
\begin{itemize}
\item Let $(X,d)$ be a pseudo-metric space. 
A singleton set $\{x\}$ is a nonstandardly closed and nonstandardly type-III-closed
set, since every ball $B(x;\epsilon )$ contains $\{x\}$, which also 
implies that $B(x;\epsilon )\oplus\Omega$
contains $\{x\}\oplus\Omega$. 

\item Let $(X,d)$ be a pseudo-metric space. 
If the singleton set $\{x\}$ has the null decomposition, then $\{x\}$ is 
also a nonstandardly type-I-closed set by Remark~\ref{ch1r396}.

\item Let $(X,d)$ be a pseudo-metric space such that the pseudo-metric $d$
satisfies the null equalities. Now, for any $\omega\in\Omega$, 
we have $d(x\oplus\omega ,x)=d(x,x)=0$,
which says that $x\oplus\omega\in B(x;\epsilon )$, i.e., the open ball
$B(x;\epsilon )$ contains points of $\{x\}\oplus\Omega$. It says that 
the singleton set $\{x\}$ is a nonstandardly type-II-closed set. 
\end{itemize}
}\end{Rem}

\begin{Pro}
Let $(X,d)$ be a pseudo-metric space on a nonstandard vector space $X$ such that the 
pseudo-metric $d$ satisfies the null inequalities. Then, we have 
\[\mbox{\em cl}(\mbox{\em cl}(A))=\mbox{\em cl}(A)\mbox{ and }
\mbox{\em cl}^{\tiny \mbox{\em (II)}}(\mbox{\em cl}^{\tiny \mbox{\em (II)}}(A))=\mbox{\em cl}^{\tiny \mbox{\em (II)}}(A).\]
If we further assume that the pseudo-metric $d$ satisfies the null equalities, 
then we have 
\[\mbox{\em cl}^{\tiny \mbox{\em (I)}}(\mbox{\em cl}^{\tiny \mbox{\em (I)}}(A))=\mbox{\em cl}^{\tiny \mbox{\em (I)}}(A)\mbox{ and }
\mbox{\em cl}^{\tiny \mbox{\em (III)}}(\mbox{\em cl}^{\tiny \mbox{\em (III)}}(A))=\mbox{\em cl}^{\tiny \mbox{\em (III)}}(A).\]
\end{Pro}
\begin{Proof}
We consider the case of nonstandardly type-I-closed set. It suffices to show 
the inclusion $\mbox{cl}^{\tiny \mbox{(I)}}(\mbox{cl}^{\tiny \mbox{(I)}}(A))\subseteq\mbox{cl}^{\tiny \mbox{(I)}}(A)$. 
Suppose that $x\in\mbox{cl}^{\tiny \mbox{(I)}}(\mbox{cl}^{\tiny \mbox{(I)}}(A))$. Then we want to claim that 
$x\in\mbox{cl}^{\tiny \mbox{(I)}}(A)$. If $x\in A$, then $x\in\mbox{cl}^{\tiny \mbox{(I)}}(A)$.
Therefore, we assume $x\not\in A$.
By definition, for any open ball $B(x;\epsilon )$, we have
\begin{equation}{\label{ch1eq264}}
\mbox{cl}^{\tiny \mbox{(I)}}(A)\cap [B(x;\epsilon )\oplus\Omega ]\neq\emptyset .
\end{equation}
Suppose that $y$ is an element in the intersection of (\ref{ch1eq264}). 
Then $y=\hat{x}\oplus\omega\in\mbox{cl}^{\tiny \mbox{(I)}}(A)$ with $\omega\in\Omega$ and
$\hat{x}\in B(x;\epsilon )$, i.e., $d(x,\hat{x})<\epsilon$. 
Since $y\in\mbox{cl}^{\tiny \mbox{(I)}}(A)$, we have $y\in A$, or $y\in A^{c}$ with 
\begin{equation}{\label{ch1*eq259}}
A\cap [B(y;\epsilon -\hat{\epsilon})\oplus\Omega ]\neq\emptyset,
\end{equation}
where $\hat{\epsilon}=d(x,\hat{x})$. Since $y\in B(x;\epsilon )\oplus\Omega$,
if $y\in A$, then we see that 
$A\cap [B(x;\epsilon )\oplus\Omega ]\neq\emptyset$, i.e., 
$x\in\mbox{cl}^{\tiny \mbox{(I)}}(A)$, since $x\not\in A$. 
Now, for $y\in A^{c}$ and any $\hat{y}\in B(y;\epsilon -\hat{\epsilon})$, 
we have $d(\hat{y},y)<\epsilon -\hat{\epsilon}$ and
\begin{align*}
d(\hat{y},x) & \leq d(\hat{y},y)+d(y,x)=d(\hat{y},y)+d(\hat{x}\oplus\omega ,x)\\
& =d(\hat{y},y)+d(\hat{x},x)\\
& =d(\hat{y},y)+\hat{\epsilon}<\epsilon -\hat{\epsilon}+\hat{\epsilon}=\epsilon ,
\end{align*}
which shows that $\hat{y}\in B(x;\epsilon )$, i.e., 
$B(y;\epsilon -\hat{\epsilon})\oplus\Omega\subseteq
B(x;\epsilon )\oplus\Omega$. From (\ref{ch1*eq259}), it also says that 
$A\cap [B(x;\epsilon )\oplus\Omega ]\neq\emptyset$, i.e.,
$x\in\mbox{cl}^{\tiny \mbox{(I)}}(A)$, since $x\not\in A$.
This shows the inclusion $\mbox{cl}^{\tiny \mbox{(I)}}(\mbox{cl}^{\tiny \mbox{(I)}}(A))\subseteq\mbox{cl}^{\tiny \mbox{(I)}}(A)$.

For the case of nonstandardly type-II-closed set, 
suppose that $x\in\mbox{cl}^{\tiny \mbox{(II)}}(\mbox{cl}^{\tiny \mbox{(II)}}(A))$. Then we want to claim that 
$x\in\mbox{cl}^{\tiny \mbox{(II)}}(A)$. As described above, we may assume that $x\not\in A$. 
By definition, for any open ball $B(x;\epsilon )$, we have
\begin{equation}{\label{ch1*eq264}}
B(x;\epsilon )\cap [\mbox{cl}^{\tiny \mbox{(II)}}(A)\oplus\Omega ]\neq\emptyset .
\end{equation}
Suppose that $y$ is an element in the intersection of (\ref{ch1*eq264}). 
Then $y=\hat{y}\oplus\omega\in B(x;\epsilon )$ with $\omega\in\Omega$ and 
$\hat{y}\in\mbox{cl}^{\tiny \mbox{(II)}}(A)$. Then we have $d(\hat{y},x)\leq
d(\hat{y}\oplus\omega ,x)=d(y,x)<\epsilon$,.
Since $\hat{y}\in\mbox{cl}^{\tiny \mbox{(II)}}(A)$, we have $\hat{y}\in A$, 
or $\hat{y}\in A^{c}$ with 
\begin{equation}{\label{ch1**eq259}}
B(\hat{y};\epsilon -\hat{\epsilon})\cap [A\oplus\Omega ]\neq\emptyset,
\end{equation}
where $\hat{\epsilon}=d(\hat{y},x)$. If $\hat{y}\in A$, then 
$y=\hat{y}\oplus\omega\in B(x;\epsilon )\cap [A\oplus\Omega ]$, i.e., 
$B(x;\epsilon )\cap [A\oplus\Omega ]\neq\emptyset$. This says that 
$x\in\mbox{cl}^{\tiny \mbox{(II)}}(A)$, since $x\not\in A$. 
Now, for $\hat{y}\in A^{c}$ and any $\hat{x}\in 
B(\hat{y};\epsilon -\hat{\epsilon})$, 
we have $d(\hat{x},\hat{y})<\epsilon -\hat{\epsilon}$ and
\[d(\hat{x},x)\leq d(\hat{x},\hat{y})+d(\hat{y},x)
<\epsilon -\hat{\epsilon}+\hat{\epsilon}=\epsilon ,\]
which shows that $\hat{x}\in B(x;\epsilon )$, i.e., 
$B(\hat{y};\epsilon -\hat{\epsilon})\subseteq B(x;\epsilon )$. 
From (\ref{ch1**eq259}), it also says that 
$B(x;\epsilon )\cap [A\oplus\Omega ]\neq\emptyset$, i.e.,
$x\in\mbox{cl}^{\tiny \mbox{(II)}}(A)$, since $x\not\in A$. Therefore, 
$x\in\mbox{cl}^{\tiny \mbox{(II)}}(A)$. This shows the inclusion $\mbox{cl}^{\tiny \mbox{(II)}}(\mbox{cl}^{\tiny \mbox{(II)}}(A))
\subseteq\mbox{cl}^{\tiny \mbox{(II)}}(A)$. Without considering the null set $\Omega$, 
the above arguments can also show 
$\mbox{cl}(\mbox{cl}(A))\subseteq\mbox{cl}(A)$. 

Finally, for the case of nonstandardly type-III-closed set,
suppose that $x\in\mbox{cl}^{\tiny \mbox{(III)}}(\mbox{cl}^{\tiny \mbox{(III)}}(A))$. Now we may assume that 
$x\not\in A$. By definition, for any open ball $B(x;\epsilon )$, we have
\begin{equation}{\label{ch1**eq264}}
[B(x;\epsilon )\oplus\Omega ]\cap [\mbox{cl}^{\tiny \mbox{(III)}}(A)\oplus\Omega ]
\neq\emptyset .
\end{equation}
Suppose that $y$ is an element in the intersection of (\ref{ch1**eq264}). 
Then $y=\hat{y}\oplus\omega_{1}=\hat{x}\oplus\omega_{2}$ with $\omega_{1},
\omega_{2}\in\Omega$, $\hat{x}\in B(x;\epsilon )$ and 
$\hat{y}\in\mbox{cl}^{\tiny \mbox{(III)}}(A)$. Since $d$ satisfies null equalities, we have 
\[d(\hat{y},x)=d(\hat{y}\oplus\omega_{1},x)=d(\hat{x}\oplus\omega_{2},x)
=d(\hat{x},x)<\epsilon ,\]
Since $\hat{y}\in\mbox{cl}^{\tiny \mbox{(III)}}(A)$, we have $\hat{y}\in A$, or 
$\hat{y}\in A^{c}$ with 
\begin{equation}{\label{ch1***eq259}}
[B(\hat{y};\epsilon -\hat{\epsilon})\oplus\Omega ]\cap 
[A\oplus\Omega ]\neq\emptyset,
\end{equation}
where $\hat{\epsilon}=d(\hat{y},x)$. If $\hat{y}\in A$, then 
$y\in [B(x;\epsilon )\oplus\Omega ]\cap 
[A\oplus\Omega ]\neq\emptyset$, i.e., $x\in\mbox{cl}^{\tiny \mbox{(III)}}(A)$, 
since $x\not\in A$. Now, for $\hat{y}\in A^{c}$ and any $z\in 
B(\hat{y};\epsilon -\hat{\epsilon})$, 
we have $d(z,\hat{y})<\epsilon -\hat{\epsilon}$ and
\[d(z,x)\leq d(z,\hat{y})+d(\hat{y},x)
<\epsilon -\hat{\epsilon}+\hat{\epsilon}=\epsilon ,\]
which shows that $z\in B(x;\epsilon )$, i.e., 
$B(\hat{y};\epsilon -\hat{\epsilon})\oplus\Omega\subseteq B(x;\epsilon )
\oplus\Omega$. From (\ref{ch1***eq259}), it also says that 
$[B(x;\epsilon )\oplus\Omega ]\cap [A\oplus\Omega ]\neq\emptyset$, i.e.,
$x\in\mbox{cl}^{\tiny \mbox{(III)}}(A)$, since $x\not\in A$. Therefore, 
$x\in\mbox{cl}^{\tiny \mbox{(III)}}(A)$. This shows the inclusion $\mbox{cl}^{\tiny \mbox{(III)}}(\mbox{cl}^{\tiny \mbox{(III)}}(A))
\subseteq\mbox{cl}^{\tiny \mbox{(III)}}(A)$. We complete the proof.
\end{Proof}

Inspired by the above proposition, we may consider the different combinations
of nonstandardly closed sets like $\mbox{cl}^{\tiny \mbox{(I)}}(\mbox{cl}^{\tiny \mbox{(II)}}(A))$, 
$\mbox{cl}^{\tiny \mbox{(III)}}(\mbox{cl}^{\tiny \mbox{(I)}}(A))$, $\mbox{cl}^{\tiny \mbox{(II)}}(\mbox{cl}^{\tiny \mbox{(III)}}(A))$ and so on.
The following results are very useful for further discussion.

\begin{Lem}{\label{ch1p601}}
Let $X$ be a nonstandard vector space over $\mathbb{F}$ and 
own the null decomposition. Suppose that $\Omega$ owns the self-decomposition. 
Let $A$ be any subset of $X$. Then the following statements hold true.
\begin{enumerate}
\item [{\em (i)}] Given any fixed $x_{0}\in X$, $x\oplus\omega\in A\oplus\Omega\oplus
x_{0}$ implies $x\in A\oplus\Omega\oplus x_{0}$ for any $\omega\in\Omega$
and $x\in X$. In particular, we have $x\oplus\omega\in A\oplus\Omega$ 
implies $x\in A\oplus\Omega$ for any $\omega\in\Omega$
and $x\in X$. In this case, we also have $A\subseteq A\oplus\Omega$.

\item [{\em (ii)}] Given any $x\in X$, we have the following properties.
\begin{itemize}
\item $x\oplus\omega\in A\oplus\omega$ implies
$x\in A\oplus\Omega$ for any $\omega\in\Omega$.

\item $x\oplus\omega\in A$ implies $x\in A\oplus\Omega$.
\end{itemize}

\item [{\em (iii)}] We further assume that $\Omega$ is closed under the vector addition.
Given any fixed $x_{0}\in X$, $x\oplus\omega\in A\oplus\Omega\oplus
x_{0}$ if and only if $x\in A\oplus\Omega\oplus x_{0}$ for any $\omega\in\Omega$
and $x\in X$. In particular, we have $x\oplus\omega\in A\oplus\Omega$ 
if and only if $x\in A\oplus\Omega$ for any $\omega\in\Omega$ and $x\in X$.
\end{enumerate}
\end{Lem}
\begin{Proof}
We have $x\oplus\omega =a\oplus\omega_{0}\oplus x_{0}$ for some $a\in A$
and $\omega_{0}\in\Omega$. By adding $-x$ on both sides, we have 
$\omega_{1}\oplus\omega =a\ominus x\oplus\omega_{0}\oplus x_{0}$, where $\omega_{1}=
x\ominus x\in\Omega$. Let $\omega_{2}=\omega_{1}\oplus\omega\in\Omega$.
Then we obtain $\omega_{2}=a\ominus x\oplus\omega_{0}\oplus x_{0}$. 
Since $X$ owns the null 
decomposition, we have $x=\hat{x}\oplus\omega_{3}$ for some $\hat{x}\in X$ 
and $\omega_{3}\in\Omega$. Then we have 
\begin{align*}
x & =\hat{x}\oplus\omega_{3}=\hat{x}\oplus\omega_{4}\oplus\omega_{2}
\mbox{ (where $\omega_{4}\in\Omega$, since $\Omega$ owns the 
self-decomposition)}\\
& =\hat{x}\oplus\omega_{4}\oplus a\ominus x\oplus\omega_{0}\oplus x_{0}=
\hat{x}\oplus a\ominus x\oplus\omega_{5}\oplus x_{0}
\mbox{ (where $\omega_{5}=\omega_{4}\oplus\omega_{0}\in\Omega$)}\\
& =\hat{x}\oplus a\ominus x\oplus\omega_{6}\oplus\omega_{3}\oplus x_{0}
\mbox{ (where $\omega_{6}\in\Omega$, since $\Omega$ owns the  
self-decomposition)}\\
& =x\oplus a\ominus x\oplus\omega_{6}\oplus x_{0}=a\oplus\omega_{7}\oplus x_{0}
\in A\oplus\Omega\oplus x_{0}
\mbox{ (where $\omega_{7}=\omega_{6}\oplus x\ominus x\in\Omega$)}.
\end{align*}
This shows that $x\oplus\omega\in A\oplus\Omega\oplus
x_{0}$ implies $x\in A\oplus\Omega\oplus x_{0}$ for any $\omega\in\Omega$
and $x\in X$.
If we further assume that $\Omega$ is closed under the vector addition, 
the converse is obvious. Without considering $x_{0}$, we can also show that 
$x\oplus\omega\in A\oplus\Omega$ implies $x\in A\oplus\Omega$. 
In this case, if $x\in A$, then $x\oplus\omega\in A\oplus\Omega$,
which implies $x\in A\oplus\Omega$, i.e., 
$A\subseteq A\oplus\Omega$. This proves (i) and (iii). 
The above arguments can also obtain the results (ii).
We complete the proof.
\end{Proof}

\begin{Pro}{\label{ch1p428}}
Let $(X,d)$ be a pseudo-metric space on a nonstandard vector space $X$ such that 
the pseudo-metric $d$ satisfies the null inequalities. 
If $X$ owns the null decomposition and $\Omega$ owns the self-decomposition
and is closed under the vector addition, then $\mbox{\em cl}^{\tiny \mbox{\em (II)}}(A)
\subseteq\mbox{\em cl}^{\tiny \mbox{\em (II)}}(A)\oplus\omega$ and $\mbox{\em cl}^{\tiny \mbox{\em (III)}}(A)
\subseteq\mbox{\em cl}^{\tiny \mbox{\em (III)}}(A)\oplus\omega$ for any $\omega\in\Omega$.
In other words, if $A$ is a nonstandardly type-II or type-III closed subset 
of $X$, then $A\subseteq A\oplus\omega$ for any $\omega\in\Omega$.
\end{Pro}
\begin{Proof}
Let $a\in\mbox{cl}^{\tiny \mbox{(II)}}(A)$. Since $a$ has the null decomposition, 
we have $a=\hat{a}\oplus\hat{\omega}_{0}$ for some $\hat{a}\in X$ and 
$\hat{\omega}_{0}\in\Omega$.
Given any $\omega\in\Omega$, since $\Omega$ owns the self-decomposition, 
we have $\hat{\omega}_{0}=\omega\oplus\hat{\omega}_{1}$ for some 
$\hat{\omega}_{1}\in\Omega$, i.e., $a=\hat{a}\oplus\hat{\omega}_{1}\oplus\omega\in A$. 
We want to claim that $\hat{a}\oplus\hat{\omega}_{1}$ is also in 
$\mbox{cl}^{\tiny \mbox{(II)}}(A)$. By definition, we just need to consider 
$\hat{a}\oplus\hat{\omega}_{1}\in A^{c}$.
Now if $a\in A$, i.e., $a=(\hat{a}\oplus\hat{\omega}_{1})\oplus\omega$, 
then $\hat{a}\oplus\hat{\omega}_{1}\in A\oplus\Omega$ by 
Lemma~\ref{ch1p601} (ii). This says that $\hat{a}\oplus\hat{\omega}_{1}
\in B(\hat{a}\oplus\hat{\omega}_{1};\epsilon )\cap (A\oplus\Omega )\neq\emptyset$, i.e., 
$\hat{a}\oplus\hat{\omega}_{1}\in\mbox{cl}^{\tiny \mbox{(II)}}(A)$. 
Suppose that $a\in A^{c}$. By definition, every open ball $B(a;\epsilon )$ 
contains points of $A\oplus\Omega$. We want to claim that 
$B(a;\epsilon )\subseteq B(\hat{a}\oplus\hat{\omega}_{1};\epsilon )$. 
For any $x\in B(a;\epsilon )$, we have $d(x,a)<\epsilon$ and 
\[d(\hat{a}\oplus\hat{\omega}_{1},x)\leq d(\hat{a}\oplus\hat{\omega}_{1},a)
+d(a,x)\leq d(\hat{a}\oplus\hat{\omega}_{1}\oplus\omega ,a)+d(a,x)
=d(a,a)+d(a,x)<0+\epsilon =\epsilon ,\]
which says that $x\in B(\hat{a}\oplus\hat{\omega}_{1};\epsilon )$. 
Therefore, we conclude that $B(\hat{a}\oplus\hat{\omega}_{1};\epsilon )$ 
contains points of $A\oplus\Omega$, i.e.,
$\hat{a}\oplus\hat{\omega}_{1}\in\mbox{cl}^{\tiny \mbox{(II)}}(A)$. This shows that 
$\mbox{cl}^{\tiny \mbox{(II)}}(A)\subseteq\mbox{cl}^{\tiny \mbox{(II)}}(A)\oplus\omega$.

Let $a\in\mbox{cl}^{\tiny \mbox{(III)}}(A)$. We want to claim that $\hat{a}\oplus\hat{\omega}_{1}$ 
is also in $\mbox{cl}^{\tiny \mbox{(III)}}(A)$. Suppose that $a\in A$.
By Proposition~\ref{ch1p273}, we see that $\hat{a}\oplus\hat{\omega}_{1}
\in (B(\hat{a}\oplus\hat{\omega}_{1};\epsilon )\oplus\Omega )
\cap (A\oplus\Omega )\neq\emptyset$, i.e., $\hat{a}\oplus\hat{\omega}_{1}
\in\mbox{cl}^{\tiny \mbox{(III)}}(A)$. Suppose that $a\in A^{c}$. By definition, 
every set $B(a;\epsilon )\oplus\Omega$ contains points of $A\oplus\Omega$. 
Since $B(a;\epsilon )\subseteq B(\hat{a}\oplus\hat{\omega}_{1};\epsilon )$, 
we see that every set $B(\hat{a}\oplus\hat{\omega}_{1};\epsilon )\oplus\Omega$ 
contains points of $A\oplus\Omega$,
i.e., $\hat{a}\oplus\hat{\omega}_{1}\in\mbox{cl}^{\tiny \mbox{(III)}}(A)$. This shows that 
$\mbox{cl}^{\tiny \mbox{(III)}}(A)\subseteq\mbox{cl}^{\tiny \mbox{(III)}}(A)\oplus\omega$.
\end{Proof}

\begin{Pro}{\label{ch1p259}}
Let $(X,d)$ be a pseudo-metric space on a nonstandard vector space $X$. 
Then the following statements hold true.
\begin{enumerate}
\item [{\em (i)}] The complement of an nonstandardly open set is nonstandardly closed, 
and the complement of a nonstandardly closed set is nonstandardly open.

\item [{\em (ii)}] The complement of a nonstandardly type-I-closed set is 
nonstandardly type-I-open, and the complement of a 
nonstandardly type-I-open set is nonstandardly type-I-closed.
\end{enumerate}
\end{Pro}
\begin{Proof}
Let $O$ be a nonstandardly type-I-open set and $O^{c}$ be its complement. 
For $x\in O=(O^{c})^{c}$, there is an open ball $B(x;\epsilon )$ such that 
$B(x;\epsilon )\oplus\Omega\subseteq O$, i.e., $B(x;\epsilon )\oplus\Omega$
is disjoint from $O^{c}$. It means that $x$ is not a nonstandard
type-I closure point of $O^{c}$. By Remark~\ref{ch1r420} (iv), we see
that $O^{c}$ is nonstandardly type-I-closed.
Similarly, without considering the null set $\Omega$, we can show that 
$O^{c}$ is a nonstandardly closed set if $O$ is an nonstandardly open set. 
On the other hand, let $A$ be a nonstandardly type-I-closed set and $A^{c}$ 
be its complement. By definition, we have $A=\mbox{cl}^{\tiny \mbox{(I)}}(A)$. 
For $x\in A^{c}$, $x$ is not a nonstandard type-I closure point of $A$. It says that 
there is an open ball $B(x;\epsilon )$ such that $B(x;\epsilon )\oplus 
\Omega$ is disjoint from $A$, i.e., 
$B(x;\epsilon )\oplus\Omega$ is contained in $A^{c}$. This shows that 
$x\in\mbox{int}^{\tiny \mbox{(I)}}(A^{c})$, i.e., 
$A^{c}\subseteq\mbox{int}^{\tiny \mbox{(I)}}(A^{c})$. 
Therefore, $A^{c}$ is nonstandardly type-I-open. Similarly, we can show 
that $A^{c}$ is an nonstandardly open set if $A$ is a nonstandardly closed set 
without considering the null set $\Omega$. 
\end{Proof}

\begin{Pro}{\label{ch2p510}}
Let $(X,d)$ be a pseudo-metric space on a nonstandard vector space $X$ such that the 
pseudo-metric $d$ satisfies the null equalities. 
Then, the following statement holds true.
\begin{enumerate}
\item [{\em (i)}] The complement of a nonstandardly type-II-open set is 
nonstandardly type-II-closed.

\item [{\em (ii)}] The complement of a nonstandardly type-III-open set is 
nonstandardly type-III-closed.
\end{enumerate}
\end{Pro}
\begin{Proof}
To prove part (i), we first claim $A^{c}\oplus\Omega\subseteq
(A\oplus\Omega )^{c}$. For any $x\in A^{c}\oplus\Omega$, we have 
$x=\hat{x}\oplus\hat{\omega}$ for some $\hat{x}\in A^{c}=(\mbox{int}^{\tiny \mbox{(II)}}(A))^{c}$
and $\hat{\omega}\in\Omega$. By definition, we see that $B(\hat{x};\epsilon )
\not\subseteq A\oplus\Omega$ for every $\epsilon >0$. Since $B(x;\epsilon )=
B(\hat{x}\oplus\hat{\omega};\epsilon )=B(\hat{x};\epsilon )$ by 
Proposition~\ref{ch2p1}, we also have $B(x;\epsilon )
\not\subseteq A\oplus\Omega$ for every $\epsilon >0$. This says that 
$x\not\in\mbox{int}(A\oplus\Omega )$. Now we want to show that 
$\mbox{int}(A\oplus\Omega )=A\oplus\Omega$. For $y\in A\oplus\Omega$,
we have $y=y_{0}\oplus\omega_{0}$ for some $y_{0}\in A=\mbox{int}^{\tiny \mbox{(II)}}(A)$
and $\omega_{0}\in\Omega$. By definition, there exists $\hat{\epsilon}>0$
such that $B(y_{0};\hat{\epsilon})\subseteq A\oplus\Omega$. By 
Proposition~\ref{ch2p1} again, we also see that 
\[B(y;\hat{\epsilon})=B(y_{0}\oplus\omega_{0};\hat{\epsilon})=
B(y_{0};\hat{\epsilon})\subseteq A\oplus\Omega ,\]
which shows that $y\in\mbox{int}(A\oplus\Omega )$, i.e., 
$\mbox{int}(A\oplus\Omega )=A\oplus\Omega$. Then we have 
$x\not\in A\oplus\Omega$, i.e., $x\in (A\oplus\Omega )^{c}$.
This shows that $A^{c}\oplus\Omega\subseteq (A\oplus\Omega )^{c}$.
Now, for $x\in A=\mbox{int}^{\tiny \mbox{(II)}}(A)$, there is an open ball $B(x;\epsilon )
\subseteq A\oplus\Omega$, i.e., the open ball is disjoint from 
$(A\oplus\Omega )^{c}$. Therefore, we conclude that there is an open ball 
$B(x;\epsilon )$ which is disjoint from $A^{c}\oplus\Omega$, i.e., 
$x$ is not a nonstandard type-II closure point of $A^{c}$. By 
Remark~\ref{ch1r420} (iv), we see
that $A^{c}$ is nonstandardly type-II-closed. 

To prove part (ii), we still need to claim the inclusion $A^{c}\oplus\Omega\subseteq
(A\oplus\Omega )^{c}$. For any $x\in A^{c}\oplus\Omega$, we have 
$x=\hat{x}\oplus\hat{\omega}$ for some $\hat{x}\in A=(\mbox{int}^{\tiny \mbox{(III)}}(A))^{c}$
and $\hat{\omega}\in\Omega$. By the arguments of (i), we see that $B(x;\epsilon )
\oplus\Omega =B(\hat{x};\epsilon )\oplus\Omega\not\subseteq A\oplus\Omega$ 
for every $\epsilon >0$. This says that $x\not\in\mbox{int}^{\tiny \mbox{(I)}}(A\oplus\Omega )$.
Now we want to show that 
$\mbox{int}^{\tiny \mbox{(I)}}(A\oplus\Omega )=A\oplus\Omega$. For $y\in A\oplus\Omega$,
we have $y=y_{0}\oplus\omega_{0}$ for some $y_{0}\in A=\mbox{int}^{\tiny \mbox{(III)}}(A)$
and $\omega_{0}\in\Omega$. By definition, there exists $\hat{\epsilon}>0$
such that $B(y;\hat{\epsilon})\oplus\Omega =
B(y_{0};\hat{\epsilon})\oplus\Omega\subseteq A\oplus\Omega$,
which shows that $y\in\mbox{int}^{\tiny \mbox{(I)}}(A\oplus\Omega )$, i.e., 
$\mbox{int}^{\tiny \mbox{(I)}}(A\oplus\Omega )=A\oplus\Omega$. Then we have 
$x\not\in A\oplus\Omega$, i.e., $x\in (A\oplus\Omega )^{c}$.
This shows that $A^{c}\oplus\Omega\subseteq (A\oplus\Omega )^{c}$.
Now, for $x\in A=\mbox{int}^{\tiny \mbox{(III)}}(A)$, there is an open ball $B(x;\epsilon )$
such that $B(x;\epsilon )\oplus\Omega\subseteq A\oplus\Omega$, i.e., 
the set $B(x;\epsilon )\oplus\Omega$ is disjoint from 
$(A\oplus\Omega )^{c}$. Therefore, we conclude that there is an open ball 
$B(x;\epsilon )$ such that the set $B(x;\epsilon )\oplus\Omega$ 
is disjoint from $A^{c}\oplus\Omega$, i.e., 
$x$ is not a nonstandard type-III closure point of $A^{c}$. By 
Remark~\ref{ch1r420} (iv), we see
that $A^{c}$ is nonstandardly type-III-closed.
\end{Proof}

\begin{Pro}{\label{ch1*p259}}
Let $X$ be a pseudo-metric space on a nonstandard vector space $X$ such that the 
pseudo-metric $d$ satisfies the null inequalities and $X$ owns the null decomposition, 
and $\Omega$ owns the self-decomposition and is closed under the vector addition. 
Then the following statements hold true.
\begin{enumerate}
\item [{\em (i)}] The complement of a nonstandardly type-II-closed set is 
simultaneously nonstandardly open and nonstandardly type-II-open.

\item [{\em (ii)}] The complement of a nonstandardly 
type-III-closed set is simultaneously nonstandardly type-I and type-III-open.
\end{enumerate}
\end{Pro}
\begin{Proof}
To prove part (i), for $x\in A^{c}$, $x$ is not a nonstandard type-II 
closure point of $A$. It says that there is an open ball 
$B(x;\epsilon )$ such that $B(x;\epsilon )$ is disjoint from 
$A\oplus\Omega$, i.e., $B(x;\epsilon )$ is contained in 
$(A\oplus\Omega )^{c}\subseteq A^{c}$ by Lemma~\ref{ch1p601} (iii). 
This says that $A^{c}$ is nonstandardly open. Since $A^{c}\subseteq A^{c}\oplus\Omega$
by Lemma~\ref{ch1p601} (iii) again, i.e, $B(x;\epsilon )$ is 
contained in $A^{c}\oplus\Omega$. This shows that 
$x\in\mbox{int}^{\tiny \mbox{(II)}}(A^{c})$, i.e., $A^{c}$ is nonstandardly type-II-open. 

To prove part (ii), for $x\in A^{c}$, $x$ is not a 
nonstandard type-III closure point of $A$. It says that there is an 
open ball $B(x;\epsilon )$ such that $B(x;\epsilon )\oplus\Omega$ is 
disjoint from $A\oplus\Omega$, i.e., $B(x;\epsilon )\oplus\Omega$ is 
contained in $(A\oplus\Omega )^{c}\subseteq A^{c}$ by 
Lemma~\ref{ch1p601} (iii). This says that $A^{c}$ is nonstandardly 
type-I open. Since $A^{c}\subseteq A^{c}\oplus\Omega$, i.e, 
$B(x;\epsilon )\oplus\Omega$ is contained in $A^{c}\oplus\Omega$. 
This shows that $x\in\mbox{int}^{\tiny \mbox{(III)}}(A^{c})$, i.e., 
$A^{c}$ is nonstandardly type-III-open. This completes the proof.
\end{Proof}

\begin{Pro}
Let $(X,d)$ be a pseudo-metric space on a nonstandard vector space $X$ such that the 
pseudo-metric $d$ satisfies the null equalities. 
If $A$ is nonstandardly closed or type-I-closed, then $A$ is simultaneously
nonstandardly closed, type-I-closed, type-II-closed and type-III-closed.
\end{Pro}
\begin{Proof}
Let $A$ be a nonstandardly closed or type-I-closed set. 
Then Proposition~\ref{ch1p259} says that 
$A^{c}$ is nonstandardly open or type-I-open, respectively,
which is also simultaneously nonstandardly open, type-I-open, type-II-open
and type-III-open by Proposition~\ref{ch1p409} (viii). Using 
Propositions~\ref{ch1p259} and \ref{ch2p510}, $A=(A^{c})^{c}$ is simultaneously
nonstandardly closed, type-I-closed, type-II-closed and type-III-closed. 
This completes the proof.
\end{Proof}

\begin{Pro}
Let $(X,d)$ be a pseudo-metric space on a nonstandard vector space $X$. 
Then the following statement holds true.
\begin{enumerate}
\item [{\em (i)}] A closed ball is also a nonstandardly closed and 
type-III-closed subset of $X$.

\item [{\em (ii)}] Suppose that the pseudo-metric $d$ satisfies the null inequalities. 
A closed ball is also a nonstandardly type-I-closed subset of $X$. 

\item [{\em (iii)}] Suppose that $X$ owns the null decomposition and 
the pseudo-metric $d$ satisfies the null equalities. 
A closed ball is also a nonstandardly type-II-closed subset of $X$.
\end{enumerate}
\end{Pro}
\begin{Proof}
To prove part (i), let $\bar{B}(x_{0};\epsilon )$ be a closed ball. 
We want to claim that its complement $\bar{B}^{c}(x_{0};\epsilon )$ is 
a nonstandardly open set.
For any $x\in\bar{B}^{c}(x_{0};\epsilon )$, we have $d(x,x_{0})>\epsilon$. 
Let $\hat{\epsilon}=d(x,x_{0})$. We consider the open ball 
$B(x,\hat{\epsilon}-\epsilon )$. Then, for any element 
$\hat{x}\in B(x;\hat{\epsilon}-\epsilon )$, i.e.,
$d(\hat{x},x)<\hat{\epsilon}-\epsilon$, we have 
\[\hat{\epsilon}=d(x,x_{0})\leq d(x,\hat{x})+d(\hat{x},x_{0})<
\hat{\epsilon}-\epsilon +d(\hat{x},x_{0}),\]
which implies $d(\hat{x},x_{0})>\epsilon$, i.e., 
$B(x,\hat{\epsilon}-\epsilon )\subseteq\bar{B}^{c}(x_{0};\epsilon )$.
This shows that $\bar{B}^{c}(x_{0};\epsilon )$ is a nonstandardly open set.
According to Proposition~\ref{ch1p259} (i), it says that
$\bar{B}(x_{0};\epsilon )$ is a nonstandardly closed set. 
We also have $B(x,\hat{\epsilon}-\epsilon )\oplus\Omega\subseteq
\bar{B}^{c}(x_{0};\epsilon )\oplus\Omega$ by adding $\Omega$ on both sides, 
i.e., $\bar{B}^{c}(x_{0};\epsilon )$ is nonstandardly type-III-open. 
Therefore, by Proposition~\ref{ch1p259} (ii), $\bar{B}(x_{0};\epsilon )$ is a 
nonstandardly type-III-closed set.  

To prove part (ii), from the above arguments, since $d(\hat{x}\oplus\omega ,x_{0})
\geq d(\hat{x},x_{0})>\epsilon$, 
we have $\hat{x}\oplus\omega\in\bar{B}^{c}(x_{0};\epsilon )$, i.e., 
$B(x,\hat{\epsilon}-\epsilon )\oplus\Omega\subseteq
\bar{B}^{c}(x_{0};\epsilon )$. This shows that 
$\bar{B}^{c}(x_{0};\epsilon )$ is nonstandardly type-I-open. Therefore,
by Proposition~\ref{ch1p259} (i), $\bar{B}(x_{0};\epsilon )$ is a 
nonstandardly type-I-closed set.

To prove part (iii), for $\hat{x}\in\bar{B}^{c}(x_{0};\epsilon )$, since $\hat{x}$ 
has the null decomposition, i.e., $\hat{x}=x_{1}\oplus\omega_{1}$ for some 
$\omega_{1}\in\Omega$, we have $d(x_{1},x_{0})=d(x_{1}\oplus\omega_{1},
x_{0})=d(\hat{x},x_{0})>\epsilon$, which says that $x_{1}\in 
\bar{B}^{c}(x_{0};\epsilon )$. Therefore, we have $\hat{x}\in
\bar{B}^{c}(x_{0};\epsilon )\oplus\Omega$, i.e., 
$\bar{B}^{c}(x_{0};\epsilon )\subseteq\bar{B}^{c}(x_{0};\epsilon )\oplus\Omega$.
Since $B(x,\hat{\epsilon}-\epsilon )\subseteq\bar{B}^{c}(x_{0};\epsilon )$
by the above arguments, we also have 
$B(x,\hat{\epsilon}-\epsilon )\subseteq\bar{B}^{c}(x_{0};\epsilon )\oplus\Omega$. 
This shows that $\bar{B}^{c}(x_{0};\epsilon )$ is nonstandardly type-II-open. By
Proposition~\ref{ch1p259} (ii), $\bar{B}(x_{0};\epsilon )$ is a 
nonstandardly type-II-closed set. This completes the proof.
\end{Proof}

\section{Topoloigcal Spaces}

Now, we are in a position to investigate the topological structure generated by 
the pseudo-metric space $(X,d)$ on a nonstandard vector space $X$.
In this case, we can obtain a nonstandard topological vector space $X$.

\begin{Thm}{\label{ch1p59}}
Let $(X,d)$ be a pseudo-metric space on a nonstandard vector space $X$. 
We denote by $\tau^{\tiny \mbox{(I)}}$ the set of all nonstandardly type-I-open
subsets of $X$. Then $(X,\tau^{\tiny \mbox{(I)}})$ is a topological space.
\end{Thm}
\begin{Proof}
By part (b) of Remark~\ref{ch1r430}, we see that $\emptyset\in\tau^{\tiny \mbox{(I)}}$ and 
$X\in\tau^{\tiny \mbox{(I)}}$. Let $A=\bigcap_{i=1}^{n}A_{i}$, where $A_{i}$ 
are nonstandardly type-I-open sets for all
$i=1,\cdots ,n$. For $x\in A$, we have $x\in A_{i}$ for all $i=1,\cdots ,n$.
Then there exist $\epsilon_{i}$ such that $B(x;\epsilon_{i})\oplus\Omega
\subseteq A_{i}$ for all $i=1,\cdots ,n$. Let $\epsilon =\min\{\epsilon_{1},
\cdots ,\epsilon_{n}\}$. Then $B(x;\epsilon )\oplus\Omega\subseteq
B(x;\epsilon_{i})\oplus\Omega\subseteq A_{i}$ for all $i=1,\cdots ,n$,
which says that $B(x;\epsilon )\oplus\Omega\subseteq\bigcap_{i=1}^{n}A_{i}
=A$. Therefore, the intersection $A$ is nonstandardly type-I-open. 
On the other hand, we 
let $A=\bigcup_{\delta}A_{\delta}$. Then $x\in A$ implies that $x\in 
A_{\delta}$ for some $\delta$. This says that $B(x;\epsilon )\oplus\Omega
\subseteq A_{\delta}\subseteq A$ for some $\epsilon >0$. Therefore, the 
union $A$ is nonstandardly type-I-open. 
\end{Proof}

\begin{Thm}
Let $(X,d)$ be a pseudo-metric space on a nonstandard vector space $X$. 
We denote by $\tau_{0}$ the set of all nonstandardly open 
subsets of $X$. Then $(X,\tau_{0})$ is a topological space.
\end{Thm}
\begin{Proof}
The empty set $\emptyset$ and $X$ are nonstandardly open by part (a) of 
Remark~\ref{ch1r430}. The remaining proof follows from the arguments of 
Theorem~\ref{ch1p59} without considering the null set $\Omega$.
\end{Proof}

Let $\tau^{\tiny \mbox{(II)}}$ and $\tau^{\tiny \mbox{(III)}}$ be the families of nonstandard type-II and 
type-III open subsets of $X$, respectively. Then $\tau^{\tiny \mbox{(II)}}$ and $\tau^{\tiny \mbox{(III)}}$
cannot be the topologies. The main reason is that the equality 
\[\left [(A_{1}\oplus\Omega )\cap (A_{2}\oplus\Omega )\right ]=
(A_{1}\cap A_{2})\oplus\Omega\] 
cannot hold true in general for any 
subsets $A_{1}$ and $A_{2}$ of $X$. However, we still have some related 
results that will be discussed below.

\begin{Lem}{\label{ch2p503}}
Let $X$ be a nonstandard vector space over $\mathbb{F}$, and 
let $A_{1}$ and $A_{2}$ be subsets of $X$. Then we have 
\[(A_{1}\cap A_{2})\oplus\Omega\subseteq\left [(A_{1}\oplus\Omega )\cap 
(A_{2}\oplus\Omega )\right ].\]
If we further assume that the null set $\Omega$ owns the self-decomposition, and 
any one of the following conditions is satisfied:
\begin{itemize}
\item $A_{2}\oplus\Omega\subseteq A_{2}$ and, for any 
$\omega\in\Omega$, $a\oplus\omega\in A_{1}\oplus\omega$ implies $a\in A_{1}$.

\item $A_{1}\oplus\Omega\subseteq A_{1}$ and, for any 
$\omega\in\Omega$, $a\oplus\omega\in A_{2}\oplus\omega$ implies 
$a\in A_{2}$. 
\end{itemize}
Then we have 
\[\left [(A_{1}\oplus\Omega )\cap (A_{2}\oplus\Omega )\right ]=
(A_{1}\cap A_{2})\oplus\Omega .\]
\end{Lem}
\begin{Proof}
For $y\in (A_{1}\cap A_{2})\oplus\Omega$, we have $y=a\oplus\omega$ with 
$a\in A_{i}$ for $i=1,2$ and $\omega\in\Omega$, which also says that 
$y\in\left [(A_{1}\oplus\Omega )\cap (A_{2}\oplus\Omega )\right ]$, i.e.,
$(A_{1}\cap A_{2})\oplus\Omega\subseteq\left [(A_{1}\oplus\Omega )\cap 
(A_{2}\oplus\Omega )\right ]$. Conversely, we assume that the first
condition is satisfied. Let $x\in (A_{1}\oplus\Omega )\cap (A_{2}\oplus\Omega )$.
Then $x=a_{1}\oplus\omega_{1}=a_{2}\oplus\omega_{2}$ for some $a_{1}\in A_{1}$,
$a_{2}\in A_{2}$ and $\omega_{1},\omega_{2}\in\Omega$. Since $\Omega$ owns
the self-decomposition, we have $\omega_{2}=\hat{\omega}_{2}\oplus
\omega_{1}$ for some $\hat{\omega}_{2}\in\Omega$. If we write 
$\hat{a}_{2}=a_{2}\oplus\hat{\omega}_{2}\in A_{2}\oplus\Omega\subseteq A_{2}$, 
then we have $\hat{a}_{2}\in A_{2}$. In this case, we have 
\[\hat{a}_{2}\oplus\omega_{1}=a_{2}\oplus\hat{\omega}_{2}\oplus\omega_{1}
=a_{2}\oplus\omega_{2}=a_{1}\oplus\omega_{1}\in A_{1}\oplus\omega_{1},\]
which shows that $\hat{a}_{2}\in A_{1}$. This says that $\hat{a}_{2}\in A_{1}\cap 
A_{2}$, i.e., $x=\hat{a}_{2}\oplus\omega_{1}\in (A_{1}\cap A_{2})\oplus\Omega$.

Now suppose that the second condition is satisfied. Since $\Omega$ owns
the self-decomposition, we have $\omega_{1}=\hat{\omega}_{1}\oplus
\omega_{2}$ for some $\hat{\omega}_{1}\in\Omega$. If we write 
$\hat{a}_{1}=a_{1}\oplus\hat{\omega}_{1}\in A_{1}\oplus\Omega\subseteq A_{1}$, 
then we have $\hat{a}_{1}\in A_{1}$. In this case, we have 
\[\hat{a}_{1}\oplus\omega_{2}=a_{1}\oplus\hat{\omega}_{1}\oplus\omega_{2}
=a_{1}\oplus\omega_{1}=a_{2}\oplus\omega_{2}\in A_{2}\oplus\omega_{2},\]
which shows that $\hat{a}_{1}\in A_{2}$. This says that $\hat{a}_{1}\in A_{1}\cap 
A_{2}$, i.e., $x=\hat{a}_{1}\oplus\omega_{2}\in (A_{1}\cap A_{2})\oplus\Omega$.
This completes the proof.
\end{Proof}

We denote by $\widetilde{\tau}^{\tiny \mbox{(II)}}$ the set of all nonstandardly type-II-open 
subsets of $X$ such that $\emptyset\in\widetilde{\tau}^{\tiny \mbox{(II)}}$ and, 
for each $\emptyset\neq A\in\widetilde{\tau}^{\tiny \mbox{(II)}}$, 
the following condition is satisfied:
\begin{equation}{\label{ch2eq506}}
\mbox{$A\oplus\Omega\subseteq A$ and, for any 
$\omega\in\Omega$, $a\oplus\omega\in A\oplus\omega$ implies $a\in A$}.
\end{equation}
We also denote by $\widetilde{\tau}^{\tiny \mbox{(III)}}$ the set of all nonstandardly 
type-III-pseudo-open subsets of $X$ such that 
$\emptyset\in\widetilde{\tau}^{\tiny \mbox{(III)}}$ and, for each 
$\emptyset\neq A\in\widetilde{\tau}^{\tiny \mbox{(III)}}$, 
the above condition (\ref{ch2eq506}) is satisfied. Of course, 
we see that $\widetilde{\tau}^{\tiny \mbox{(II)}}\subseteq\tau^{\tiny \mbox{(II)}}$ and 
$\widetilde{\tau}^{\tiny \mbox{(III)}}\subseteq\tau^{\tiny \mbox{(III)}}$.

\begin{Pro}{\label{ch1p267*}}
Let $(X,d)$ be a pseudo-metric space on a nonstandard vector space $X$ such that the 
pseudo-metric $d$ satisfies the null inequalities, and let $B(x;\epsilon )$ be 
any open ball centered at $x$ with radius $\epsilon$. 
Then, the following statements hold true.
\begin{enumerate}
\item [{\em (i)}] Suppose that $X$ owns the null decomposition and the 
pseudo-metric $d$ satisfies the null equalities. Then $B(x_{0};\epsilon )
\in\widetilde{\tau}^{\tiny \mbox{(II)}}$.

\item [{\em (ii)}] Suppose that the pseudo-metric $d$ satisfies the null equalities.
Then $B(x_{0};\epsilon )\in\widetilde{\tau}^{\tiny \mbox{(III)}}$.
\end{enumerate}
\end{Pro}
\begin{Proof}
To prove part (i), from Proposition~\ref{ch1p267} (iii), we remain to show that 
$B(x;\epsilon )$ satisfies condition (\ref{ch2eq506}). From Proposition~\ref{ch2p91},
we immediately have $B(x;\epsilon )\oplus\Omega\subseteq B(x;\epsilon )$.
Now, for $a\oplus\omega\in B(x;\epsilon )\oplus\omega$, we have 
$a\oplus\omega =\hat{a}\oplus\omega$ for some $\hat{a}\in B(x;\epsilon )$,
i.e., $d(\hat{a},x)<\epsilon$. Since $d$ satisfies the null equalities, we have 
\[d(a,x)=d(a\oplus\omega ,x)=d(\hat{a}\oplus\omega ,x)=d(\hat{a},x)<\epsilon ,\]
which says that $a\in B(x;\epsilon )$. Therefore, we obtain the desired result.

To prove part (ii), from Proposition~\ref{ch1p267} (i), we remain to show that 
$B(x;\epsilon )$ satisfies condition (\ref{ch2eq506}). The arguments of (i) are also 
valid to show that $B(x;\epsilon )$ satisfies condition (\ref{ch2eq506}).
We complete the proof.
\end{Proof}

\begin{Thm}{\label{ch2p504}}
Let $(X,d)$ be a pseudo-metric space on a nonstandard vector space $X$.
Then $(X,\widetilde{\tau}^{\tiny \mbox{(II)}})$ is a topological space.
\end{Thm}
\begin{Proof}
Given $A_{1},A_{2}\in\widetilde{\tau}^{\tiny \mbox{(II)}}$, we let $A=A_{1}\cap A_{2}$. 
For $x\in A$, we have $x\in A_{i}$ for $i=1,2$.
Then there exist $\epsilon_{i}$ such that $B(x;\epsilon_{i})\subseteq 
A_{i}\oplus\Omega$ for all $i=1,2$. 
Let $\epsilon =\min\{\epsilon_{1},\epsilon_{2}\}$. 
Then $B(x;\epsilon )\subseteq B(x;\epsilon_{i})\subseteq A_{i}\oplus\Omega$ 
for all $i=1,2$, which says that $B(x;\epsilon )\subseteq\left [(A_{1}
\oplus\Omega )\cap (A_{2}\oplus\Omega )\right ]=(A_{1}\cap A_{2})
\oplus\Omega =A\oplus\Omega$ by Lemma~\ref{ch2p503}. This shows that 
$A$ is nonstandardly type-II-open. We also need to show that $A$ satisfies
condition (\ref{ch2eq506}). For $x\in A\oplus\Omega$, we have $x=a\oplus\omega$
for some $a\in A$ and $\omega\in\Omega$. Since $a\in A_{1}\cap A_{2}$,
it also says that $x\in A_{1}\oplus\Omega\subseteq A_{1}$ and 
$x\in A_{2}\oplus\Omega\subseteq A_{2}$. Therefore, we have $x\in
A_{1}\cap A_{2}=A$, i.e., $A\oplus\Omega\subseteq A$. On the other hand,
$a\oplus\omega\in A\oplus\omega =(A_{1}\cap A_{2})\oplus\omega\subseteq
A_{1}\oplus\omega$, which says that $a\in A_{1}$. We can similarly have 
$a\oplus\omega\subseteq A_{2}\oplus\omega$, which also implies $a\in A_{2}$.
Therefore, we have $a\in A_{1}\cap A_{2}=A$. This says that $A$ is indeed
in $\widetilde{\tau}^{\tiny \mbox{(II)}}$. Therefore, the intersection of finitely many 
members of $\widetilde{\tau}^{\tiny \mbox{(II)}}$ is a member of $\widetilde{\tau}^{\tiny \mbox{(II)}}$. 
Now given $\{A_{\delta}\}\subset\widetilde{\tau}^{\tiny \mbox{(II)}}$, we 
let $A=\bigcup_{\delta}A_{\delta}$. Then $x\in A$ implies that $x\in 
A_{\delta}$ for some $\delta$. This says that $B(x;\epsilon )\subseteq 
A_{\delta}\oplus\Omega\subseteq A\oplus\Omega$ for some $\epsilon >0$. 
Therefore, the union $A$ is nonstandardly type-II-open. We also need to 
show that $A$ satisfies condition (\ref{ch2eq506}). For $x\in A\oplus\Omega$,
we have $x=a\oplus\omega$, where $a\in A$, i.e., $a\in A_{\delta}$ for some
$\delta$. It also says that $x\in A_{\delta}\oplus\Omega\subseteq A_{\delta}
\subseteq A$, i.e., $A\oplus\Omega\subseteq A$. On the other hand, for 
$a\oplus\omega\in A\oplus\omega$, we have $a\oplus\omega =\hat{a}\oplus\omega$
for some $\hat{a}\in A$, i.e., $\hat{a}\in A_{\delta}$ for some $\delta$,
which also implies $a\oplus\omega =\hat{a}\oplus\omega\in A_{\delta}\oplus\omega$.
Therefore, we obtain $a\in A_{\delta}\subseteq A$. This shows that $A$ is indeed
in $\widetilde{\tau}^{\tiny \mbox{(II)}}$.
By part (c) of Remark~\ref{ch1r430}, we see that $\emptyset$ and $X$ are also 
nonstandardly type-II open subsets of $X$. It is not hard to see that 
$X$ satisfies condition (\ref{ch2eq506}). This shows that 
$X\in\widetilde{\tau}^{\tiny \mbox{(II)}}$. We complete the proof.
\end{Proof}

\begin{Thm}
Let $(X,d)$ be a pseudo-metric space on a nonstandard vector space $X$.
Then $(X,\widetilde{\tau}^{\tiny \mbox{(III)}})$ is a topological space.
\end{Thm}
\begin{Proof}
By part (d) of Remark~\ref{ch1r430}, we see that $\emptyset,X\in
\widetilde{\tau}^{\tiny \mbox{(III)}}$, since $X$ satisfies condition (\ref{ch2eq506}). 
Given $A_{1},A_{2}\in\widetilde{\tau}^{\tiny \mbox{(III)}}$, we let $A=A_{1}\cap A_{2}$. 
For $x\in A$, there exist 
$\epsilon_{i}$ such that $B(x;\epsilon_{i})\oplus\Omega\subseteq 
A_{i}\oplus\Omega$ for all $i=1,2$. Let $\epsilon =\min\{\epsilon_{1},
\epsilon_{2}\}$. Then $B(x;\epsilon )\oplus\Omega\subseteq 
B(x;\epsilon_{i})\oplus\Omega\subseteq A_{i}\oplus\Omega$ 
for all $i=1,2$, which says that $B(x;\epsilon )\oplus\Omega\subseteq
\left [(A_{1}\oplus\Omega )\cap (A_{2}\oplus\Omega )\right ]=
(A_{1}\cap A_{2})\oplus\Omega =A\oplus\Omega$ by Lemma~\ref{ch2p503}. 
This shows that $A$ is nonstandardly type-III open. From the arguments of 
Theorem~\ref{ch2p504}, we see that $A$ satisfies condition (\ref{ch2eq506}).
Therefore, the intersection of finitely many members of 
$\widetilde{\tau}^{\tiny \mbox{(III)}}$ is a member of $\widetilde{\tau}^{\tiny \mbox{(III)}}$. 
Now, given $\{A_{\delta}\}\subset\widetilde{\tau}^{\tiny \mbox{(III)}}$, we 
let $A=\bigcup_{\delta}A_{\delta}$. Then $x\in A$ implies that $x\in 
A_{\delta}$ for some $\delta$. This says that $B(x;\epsilon )\oplus\Omega
\subseteq A_{\delta}\oplus\Omega\subseteq A\oplus\Omega$ for some $\epsilon >0$. 
Therefore, the union $A$ is nonstandardly type-III open. From the arguments of 
Theorem~\ref{ch2p504}, we also see that $A$ satisfies condition 
(\ref{ch2eq506}). We complete the proof.
\end{Proof}

\begin{Rem}
{\em 
Let $(X,d)$ be a pseudo-metric space on a nonstandard vector space $X$ such that the 
pseudo-metric $d$ satisfies the null equalities and $\Omega$ 
owns the self-decomposition. Then we have the following observations.
\begin{itemize}
\item From part (ii) of Proposition~\ref{ch1p409}, 
we see that if $A$ is $\tau_{0}$-open,
then $A$ is also $\tau^{\tiny \mbox{(I)}}$-open, and if $A$ is 
$\widetilde{\tau}^{\tiny \mbox{(II)}}$-open, then $A$ is also $\widetilde{\tau}^{\tiny \mbox{(III)}}$-open. 
This says that $\tau^{\tiny \mbox{(I)}}$ is finer than $\tau_{0}$ and 
$\widetilde{\tau}^{\tiny \mbox{(III)}}$ is finer than $\widetilde{\tau}^{\tiny \mbox{(II)}}$. 

\item Suppose that $X$ owns the null decomposition. 
Then part (vii) of Proposition~\ref{ch1p409} says that $\tau_{0}=\tau^{\tiny \mbox{(I)}}$,
$\widetilde{\tau}^{\tiny \mbox{(II)}}=\widetilde{\tau}^{\tiny \mbox{(III)}}$ and $\tau^{\tiny \mbox{(II)}}=\tau^{\tiny \mbox{(III)}}$.
\end{itemize}
}\end{Rem} 

Let $(X,d)$ be a pseudo-metric space on a nonstandard vector space $X$ such that the 
pseudo-metric $d$ satisfies the null inequalities.
We denote by $\mbox{p}\tau^{\tiny \mbox{(II)}}$ the set of all nonstandardly type-II-pseudo-open 
subsets of $X$ and by $\mbox{p}\tau^{\tiny \mbox{(III)}}$ the set of all nonstandardly 
type-III-pseudo-open subsets of $X$.

\begin{Lem}{\label{ch2p3}}
Let $(X,d)$ be a pseudo-metric space on a nonstandard vector space $X$ such that the 
pseudo-metric $d$ satisfies the null inequalities.
\begin{enumerate}
\item [{\em (i)}] Suppose that $A_{1},A_{2}\in\mbox{\em p}\tau^{\tiny \mbox{\em (II)}}$. Then, we have 
\[(A_{1}\cap A_{2})\oplus\Omega\subseteq\left [(A_{1}\oplus\Omega )\cap 
(A_{2}\oplus\Omega )\right ].\] 
If we further assume that $\Omega$ owns the self-decomposition and the 
pseudo-metric $d$ satisfies the null equalities, then 
\[\left [(A_{1}\oplus\Omega )\cap (A_{2}\oplus\Omega )\right ]=
(A_{1}\cap A_{2})\oplus\Omega .\]

\item [{\em (ii)}] Suppose that $A_{1},A_{2}\in\mbox{\em p}\tau^{\tiny \mbox{\em (III)}}$.
Then we have the same results as given in {\em (i)}.
\end{enumerate}
\end{Lem}
\begin{Proof}
The results follow immediately from Proposition~\ref{ch2p503} and part (ii) of 
Proposition~\ref{ch2p2}.
\end{Proof}

\begin{Thm}
Let $(X,d)$ be a pseudo-metric space on a nonstandard vector space $X$ such that the 
pseudo-metric $d$ satisfies the null equalities and
$\Omega$ owns the self-decomposition. Then $(X,\mbox{\em p}\tau^{\tiny \mbox{\em (II)}})$ 
is a topological space.
\end{Thm}
\begin{Proof}
Given $A_{1},A_{2}\in\mbox{p}\tau^{\tiny \mbox{(II)}}$, we let $A=A_{1}\cap A_{2}$. 
We want to show $A=\mbox{pint}^{\tiny \mbox{(II)}}(A)$.
For $x\in A$, we have $x\in A_{i}$ for $i=1,2$.
Then there exist $\epsilon_{i}$ such that $B(x;\epsilon_{i})\subseteq 
A_{i}\oplus\Omega$ for all $i=1,2$. 
Let $\epsilon =\min\{\epsilon_{1},\epsilon_{2}\}$. 
Then $B(x;\epsilon )\subseteq B(x;\epsilon_{i})\subseteq A_{i}\oplus\Omega$ 
for all $i=1,2$, which says that $B(x;\epsilon )\subseteq\left [(A_{1}
\oplus\Omega )\cap (A_{2}\oplus\Omega )\right ]=(A_{1}\cap A_{2})
\oplus\Omega =A\oplus\Omega$ by Lemma~\ref{ch2p3}. This shows that 
$x\in\mbox{pint}^{\tiny \mbox{(II)}}(A)$, i.e., $A\subseteq\mbox{pint}^{\tiny \mbox{(II)}}(A)$. 
Conversely, for $x\in\mbox{pint}^{\tiny \mbox{(II)}}(A)$, 
by part (ii) of Proposition~\ref{ch2p2}, we have 
\[x\in B(x;\epsilon )\subseteq A\oplus\Omega =(A_{1}\cap A_{2})\oplus\Omega
\subseteq A_{1}\oplus\Omega\subseteq A_{1}.\]
We can similary obtain $x\in A_{2}$, i.e., $x\in A_{1}\cap A_{2}=A$.
This shows that $\mbox{pint}^{\tiny \mbox{(II)}}(A)\subseteq A$.
Therefore, the intersection of finitely many members of 
$\mbox{p}\tau^{\tiny \mbox{(II)}}$ is a member of $\mbox{p}\tau^{\tiny \mbox{(II)}}$. 
Now, given $\{A_{\delta}\}\subset\mbox{p}\tau^{\tiny \mbox{(II)}}$, we 
let $A=\bigcup_{\delta}A_{\delta}$. Then $x\in A$ implies that $x\in 
A_{\delta}$ for some $\delta$. This says that $B(x;\epsilon )\subseteq 
A_{\delta}\oplus\Omega\subseteq A\oplus\Omega$ for some $\epsilon >0$. 
Therefore, we obtain $A\subseteq\mbox{pint}^{\tiny \mbox{(II)}}(A)$. Conversely,
for $x\in\mbox{pint}^{\tiny \mbox{(II)}}(A)$, we have $x\in B(x;\epsilon )\subseteq A\oplus\Omega$,
i.e., $x=a\oplus\omega$ for some $a\in A$ and $\omega\in\Omega$.
Then we have $a\in A_{\delta}$ for some $\delta$, i.e., $x=a\oplus\omega
\in A_{\delta}\oplus\Omega\subseteq A_{\delta}\subseteq A$ 
by part (ii) of Proposition~\ref{ch2p2}. This shows that $\mbox{pint}^{\tiny \mbox{(II)}}(A)
\subseteq A$, i.e., the union $A$ is indeed in $\mbox{pint}^{\tiny \mbox{(II)}}(A)$.
By part (c) of Remark~\ref{ch1r430}, we also see that 
$\emptyset, X\in\mbox{p}\tau^{\tiny \mbox{(II)}}$. This completes the proof.
\end{Proof}

\begin{Thm}
Let $(X,d)$ be a pseudo-metric space on a nonstandard vector space $X$ such that the 
pseudo-metric $d$ satisfies the null inequalities, $X$ owns the null 
decomposition and $\Omega$ owns the self-decomposition.
Then $(X,\mbox{\em p}\tau^{\tiny \mbox{\em (III)}})$ is a topological space.
\end{Thm}
\begin{Proof}
By part (d) of Remark~\ref{ch1r430}, we also see that $\emptyset ,X\in
\mbox{p}\tau^{\tiny \mbox{(III)}}$. Given $A_{1},A_{2}\in\mbox{p}\tau^{\tiny \mbox{(III)}}$, 
we let $A=A_{1}\cap A_{2}$. For $x\in A$, there exist 
$\epsilon_{i}$ such that $B(x;\epsilon_{i})\oplus\Omega\subseteq 
A_{i}\oplus\Omega$ for all $i=1,2$. Let $\epsilon =\min\{\epsilon_{1},
\epsilon_{2}\}$. Then $B(x;\epsilon )\oplus\Omega\subseteq 
B(x;\epsilon_{i})\oplus\Omega\subseteq A_{i}\oplus\Omega$ 
for all $i=1,2$, which says that $B(x;\epsilon )\oplus\Omega\subseteq
\left [(A_{1}\oplus\Omega )\cap (A_{2}\oplus\Omega )\right ]=
(A_{1}\cap A_{2})\oplus\Omega =A\oplus\Omega$ by Lemma~\ref{ch2p3}. 
This shows that $x\in\mbox{p}\tau^{\tiny \mbox{(III)}}$, i.e., $A\subseteq\mbox{p}\tau^{\tiny \mbox{(III)}}$.
Conversely, for $x\in\mbox{pint}^{\tiny \mbox{(III)}}(A)$, 
by Proposition~\ref{ch1p273} and part (ii) of Proposition~\ref{ch2p2}, we have 
\[x\in B(x;\epsilon )\subseteq B(x;\epsilon )\oplus\Omega\subseteq 
A\oplus\Omega =(A_{1}\cap A_{2})\oplus\Omega
\subseteq A_{1}\oplus\Omega\subseteq A_{1}.\]
We can similary obtain $x\in A_{2}$, i.e., $x\in A_{1}\cap A_{2}=A$.
This shows that $\mbox{pint}^{\tiny \mbox{(II)}}(A)\subseteq A$.
Therefore, the intersection of finitely many members of 
$\mbox{p}\tau^{\tiny \mbox{(III)}}$ is a member of $\mbox{p}\tau^{\tiny \mbox{(III)}}$. 
Now, given $\{A_{\delta}\}\subset\mbox{p}\tau^{\tiny \mbox{(III)}}$, we 
let $A=\bigcup_{\delta}A_{\delta}$. Then $x\in A$ implies that $x\in 
A_{\delta}$ for some $\delta$. This says that $B(x;\epsilon )\oplus\Omega
\subseteq A_{\delta}\oplus\Omega\subseteq A\oplus\Omega$ for some $\epsilon >0$. 
Therefore, we obtain $A\subseteq\mbox{pint}^{\tiny \mbox{(III)}}(A)$. Conversely,
for $x\in\mbox{pint}^{\tiny \mbox{(III)}}(A)$, we have $x\in B(x;\epsilon )\subseteq
B(x;\epsilon )\oplus\Omega\subseteq A\oplus\Omega$,
i.e., $x=a\oplus\omega$ for some $a\in A$ and $\omega\in\Omega$.
Then we have $a\in A_{\delta}$ for some $\delta$, i.e., $x=a\oplus\omega
\in A_{\delta}\oplus\Omega\subseteq A_{\delta}\subseteq A$ 
by part (ii) of Proposition~\ref{ch2p2}. This shows that $\mbox{pint}^{\tiny \mbox{(III)}}(A)
\subseteq A$, i.e., the union $A$ is indeed in $\mbox{pint}^{\tiny \mbox{(III)}}(A)$.
We complete the proof.
\end{Proof}

\end{document}